\documentclass[11pt]{amsart}
\usepackage{amsmath, amssymb}
 \usepackage{mathrsfs}
 
 \usepackage{ucs}
\usepackage[utf8x]{inputenc}
\usepackage{amsmath}
\usepackage{amsfonts}
\usepackage[USenglish]{babel}
\usepackage{fontenc}
\usepackage{graphicx}
\usepackage{caption}
\usepackage{subcaption}
\usepackage{amsthm}

\newcommand{\no}[1]{#1}
\renewcommand{\no}[1]{}
\no{\usepackage{times}
\usepackage[subscriptcorrection,
 slantedGreek, nofontinfo]{mtpro}
\renewcommand{\Delta}{\upDelta}}
\usepackage{color}

\date{\today}
\newtheorem{proposition}{Proposition}[section]
\newtheorem{theorem}{Theorem}[section]

\newcommand{\be}{\begin{equation}}
\newcommand{\ee}{\end{equation}}
\newcommand{\ba}{\begin{array}}
\newcommand{\ea}{\end{array}}

\newcommand{\bea}{\begin{eqnarray*}}
\newcommand{\eea}{\end{eqnarray*}}
\newcommand{\bean}{\begin{eqnarray}}
\newcommand{\eean}{\end{eqnarray}}

\def\tilde{\widetilde}

\def\cydot{\leavevmode\raise.4ex\hbox{.}}

\begin{document}
\title[multifrequency electric impedance tomography]{Numerical 
determination of anomalies in multifrequency electrical impedance tomography}

\author[Habib Ammari]{Habib Ammari$^\P$ }

%\author{}
\address{Department of Mathematics, ETH Z\"urich, 
R\"amistrasse 101, 8092 Z\"urich, Switzerland}

\email{habib.ammari@math.ethz.ch}

\author[Faouzi Triki]{Faouzi Triki$^\dag$}
\address{Laboratoire Jean Kuntzmann,  UMR CNRS 5224, 
Universit\'e  Grenoble-Alpes, 700 Avenue Centrale,
38401 Saint-Martin-d'H\`eres, France}

\email{faouzi.triki@univ-grenoble-alpes.fr}

\author[Chun-Hsiang Tsou]{ Chun-Hsiang Tsou$^\ddag$}
\address{Laboratoire Jean Kuntzmann,  UMR CNRS 5224, 
Universit\'e  Grenoble-Alpes, 700 Avenue Centrale,
38401 Saint-Martin-d'H\`eres, France}

\email{chun-hsiang.tsou@univ-grenoble-alpes.fr}

\date{\today}

\maketitle 
\begin{abstract}
The multifrequency electrical impedance tomography consists in retrieving 
 the conductivity distribution of a sample by injecting a finite number of currents 
 with multiple frequencies. In this paper we consider the case where the conductivity 
 distribution is piecewise constant, takes a constant value outside a single 
 smooth anomaly, and a frequency dependent function inside the anomaly itself.  
  Using an original spectral 
 decomposition of the solution of the forward conductivity problem in terms of
 Poincar\'e variational eigenelements,  we retrieve the Cauchy 
 data corresponding to the extreme case of a perfect conductor, and the conductivity profile. 
 We then reconstruct the anomaly from the Cauchy data. The numerical experiments
 are conducted using gradient  descent optimization algorithms.   
 
\end{abstract}

%\tableofcontents
 
 %%%%%%%%%%%%%%%%%%%%%%%%%%%%%%%%%%%%%%%%%%%%
 %%%%%%%%%%%%%%%%%%%%%%%%%%%%%%%%%%%%%%%%%%%%
\section{The mfEIT Mathematical Model}
%%%%%%%%%%%%%%%%%%%%%%%%%%%%%%%%%%%%

Experimental research has found that the conductivity of many biological tissues 
varies strongly with respect to the frequency of the applied  electric current  within certain 
frequency ranges \cite{GPG}.  In \cite{AGGJS}, using homogenization techniques, the authors 
analytically exhibit the fundamental mechanisms underlying the fact that effective biological
 tissue electrical 
properties and their frequency dependence reflect the tissue composition and physiology.
The multifrequency electrical impedance tomography (mfEIT) is a diffusive imaging 
modality that recovers the conductivity distribution of the 
tissue by using electrodes to measure the resulting voltage on
 its boundary, induced by  two known injected currents and for many frequency 
 values. The principal idea behind the  (mfEIT)  is that  the dependance of the 
effective  conductivity of the tissue with respect to the frequency of the 
 electric current  is extremely related to its state. In fact, its frequency dependence changes
 with its composition, membrane characteristics, intra-and extra-cellular fluids and 
 other factors \cite{AGGJS}. Therefore, the  frequency dependence of the conductivity of the tissue
 can  provide some information about the tissue microscopic structure and its
  physiological and pathological conditions. In other words, the frequency dependence of the 
  conductivity of the tissue can  help to determine if it is  healthy or cancerous.  The advantages
  of  the (mfEIT)  is  canceling out errors due to boundary shape, the electrode positions, and 
  other systematic errors that appear in -the more conventional imaging modality- electric
  impedance tomography (EIT) \cite{Bor}.   \\ 
In the following  we introduce the mathematical model of the (mfEIT).
Let $\Omega$ be the open bounded  smooth domain in $\mathbb R^2$, 
occupied by the sample under investigation and
denote by $\partial \Omega$ its boundary. The mfEIT forward problem 
is to determine the potential $u(\cdot, \omega) \in H^1(\Omega)
:=\{ v \in L^2(\Omega): \nabla v \in L^2(\Omega)\}$,
solution to 
\bean \label{maineq}
\left \{
\ba{lllcc}
-\nabla \cdot \left(\sigma(x, \omega) \nabla u(x, \omega)\right) =  0 & \textrm{in}
\quad \Omega, \\
\sigma(x, \omega) \partial_{\nu_\Omega} u(x, \omega)(x) = f(x) & \textrm{on}
\quad \partial \Omega,\\
\int_{\partial \Omega} u(x, \omega) ds = 0,&
\ea
\right.
\eean
where $\omega$ denotes the frequency, $\nu_\Omega(x)$ is the outward 
normal vector
to $\partial \Omega$,  $\sigma(x, \omega)$ is 
the conductivity distribution,
and   $f\in H^{-\frac{1}{2}}_\diamond(\partial \Omega) :=\{ g \in 
 H^{-\frac{1}{2}}(\partial \Omega): \int_{\partial \Omega} g \, ds =0\}$  is the
  input current. \\

In this work we are interested in the case where  the frequency dependent 
conductivity distribution takes the form 
\bean \label{conductivitydistribution}
\sigma(x, \omega) = k_0+ (k(\omega)-k_0) \chi_D(x)
\eean
with $\chi_D(x)$ being the characteristic function of a $C^2$ domain
$D$ in $\Omega$ ($\overline D\subset \Omega$), $k_0$ being a fixed 
strictly positive 
constant, and 
{$k(\omega): \mathbb
R_+\rightarrow
\mathbb C\setminus \overline{\mathbb R_-}$,} being a  continuous 
complex-valued function. \\

Here  $k_0$  represents 
 the conductivity of the background medium, is known, and 
$k(\omega)$  is the  conductivity of the
biological tissue,  given by the empirical model  
\bean \label{empirical}
k(\omega) :=  \kappa_1-\frac{\kappa_2}{ \omega^2+i\omega \kappa_3},
\eean
where  $\kappa_p>0, \; p=1, 2, 3,$ are  constants that only
depend on the biological tissue properties (see for instance \cite{AGGJS} ). 
The frequency profile $k(\omega)$ is somehow a meromorphic approximation  
with a single pole of  the graph of experimental measurements  for a  given biological
 tissue \cite{AGGJS}. It also appears as a homogenized model for periodically 
 distributed biological  cells in the dilute limit \cite{AGGJS},
 and is similar to Drude models that  describes the frequency dependence of the electric
 permittivity of a real metal within the visible frequency range \cite{MFZ}. \\

The  mfEIT inverse problem is to determine 
 the anomaly  $D$  and the characteristics $\kappa_p, \; p=1, 2, 3,$ 
 of the biological tissue   
from measurements of the boundary voltages 
$u(x, \omega)$ on  $
\partial \Omega $, for $\omega \in (\underline \omega, \overline \omega)$, 
$0\leq \underline \omega < \overline \omega$. \\

There have been  several numerical approaches on multifrequency electrical impedance 
tomography. Most of them are dealing  with small-volume anomalies \cite{ABG, ABGW, GH} and frequency-difference imaging \cite{JS, MSHA}. In small-volume-volume imaging, only the location and the multi-frequency polarization tensor can be reconstructed from boundary measurements. In frequency-difference imaging,  the main idea
is to compare  the images for different frequencies, and  consider only  the frequency
dependent part.  It was numerically shown that the approach 
can accommodate geometrical errors, including imperfectly known boundary.  This approach 
which seems more natural since it aims to identify the anomaly by  only focusing  on the
changes in the resulting images  for different frequencies, is somehow  simultaneously 
complementary and opposite  to our analysis in this paper. Taking the difference between
two images associated to  different frequencies will remove the frequency independent part
which is the keystone of the identification path pursued   in this paper.   Our strategy is based
on the plasmonic spectral decomposition  derived in \cite{AT} which splits the electric potential on the
boundary $\partial \Omega$ into two parts $u = k_0^{-1} u_0+u_f$, and  separate between the frequency 
dependent  and independent parts; see also \cite{ADM,AMRZ}. In fact the part $k_0^{-1} u_0$ corresponds   to the response of
the same  anomaly  filled with a perfect conductor, that is $k_0^{-1} u_0$ is the limit of $u$ when
$k$ tends to infinity. Precisely, in \cite{AT}, it was proven that the convergence of $u$ to  $k_0^{-1} u_0$ 
is linear in $1/k$.  We first  process algebraically the data on the boundary  and recover the frequency dependent part in order to  acquire the Cauchy data  of    frequency independent
part $u_0$.  Then, based on known results and approaches  on the reconstruction of  zero level set of harmonic  functions from Cauchy data we  determine the anomaly itself.  The other approach is
based on  perturbation techniques, suppose that the contrast $k(\omega)/k_0$ is close to
$1$, and linearize the problem around the harmonic function in the whole domain
that shares the same flux on the boundary \cite{AAJS}.  Then, the  sparsity issue comes
 into play  and  help to  speed up the convergence of the  iterative algorithm.  In our approach, the 
 notion of sparsity  appears naturally in post-processing the data on the boundary  to recover the frequency
 independent part. Precisely, it seems that only a finite number of eigenfunctions intervene in the linear inversion, and  this can be completely  mastered by the shape of the anomaly and  its distance to the boundary where the measurement are taken.  Finally, the perturbation  approach is also complementary to our analysis since the higher the frequency is the better the recovery of the frequency independent  part
 on the boundary is.  \\
 
 The paper is organized as follows.  In section 2 we provide  the spectral decomposition 
 derived in \cite{AT}.  The linearization of the frequency independent part with respect to the
  shape of the anomaly which is necessary in our identification approach  is studied in section 3.  
   Section 4 is devoted to the retrieval of the frequency independent part
 on the boundary.  Here, we will not follow the theoretical approach developed in \cite{AT} based on the unique  continuation of meromorphic functions. We solve the problem using algebraic tools under 
 simplification assumptions inspired by the sparsity properties of the conductivity distribution  and 
 the behavior of the eigenvalues near the unique accumulation point $1/2$.  We also  determine
in the sequel the profile constants $\kappa_i, \, i=1, 2, 3$. In section 5, one we
 have the Cauchy data of the frequency independent part $u_0$  on the measurement boundary 
 we use a conventional  optimization technique to recover the anomaly. Several numerical 
 examples are presented in section 6. Comments on the obtained results  and future directions
are given in the conclusion section 7.
 
%%%%%%%%%%%%%%%%%%%%%%%%%%%%%%%
\section{Spectral decomposition of  $u(x, \omega)$}
%%%%%%%%%%%%%%%%%%%%%%%%%%%%%%%
We first introduce an operator whose  spectral decomposition 
will be later  the corner stone  of the identification of the anomaly $D$.
Let $H^1_\diamond(\Omega)$ be the space of functions $v$ in  
$ H^1(\Omega)$ satisfying $\int_{\partial \Omega } v ds = 0$.\\

For $u\in H^1_\diamond(\Omega)$, we infer from the Riesz theorem that 
there exists 
a unique function $Tu\in H^1_\diamond(\Omega)$ such that for 
all $v\in H^1_\diamond(\Omega)$,  
\bean
\int_\Omega \nabla Tu \cdot \nabla v dx = \int_D \nabla u \cdot \nabla v dx. 
\eean
The variational Poincar\'e operator $T: H^1_\diamond(\Omega)\rightarrow
 H^1_\diamond(\Omega)$
is easily seen to be self-adjoint 
and bounded with norm $\|T\| \leq 1$. \\

The spectral problem for $T$  reads
as: 
Find $(\lambda, w) \in \mathbb R\times H^1_\diamond(\Omega)$, $w\not=0$
such that $\forall v \in H^1_\diamond(\Omega)$,
\bea
\lambda \int_\Omega \nabla w \cdot \nabla v dx = \int_D \nabla w \cdot \nabla v dx.
\eea

Integrating by parts, one immediately obtains that any eigenfunction 
$w$ is harmonic in $D$ and in $D^\prime = \Omega\setminus \overline D$, and
satisfies the transmission  and boundary conditions
\bea
w|^+_{\partial D} = w|^-_{\partial D}, \qquad 
\partial_{\nu_D} w|^+_{\partial D}  = 
(1-\frac{1}{\lambda})\partial_{\nu_D} w|^-_{\partial D}, \qquad  
 \partial_{\nu_\Omega} w = 0,
\eea
where $w|^\pm_{\partial D} (x)= 
\lim_{t\rightarrow 0} w(x\pm t\nu_D(x))$ for $x\in\partial D$.   
In other words, $w$ is a solution to \eqref{maineq} for 
$k = k_0(1-\frac{1}{\lambda})$ and $f=0$.\\

Let $\mathfrak H_\diamond$  the space  of harmonic functions in $D$ and
$D^\prime$, with zero mean $\int_{\partial \Omega} u ds(x) =0$, and
zero normal derivative $\partial_{\nu_\Omega} u = 0$ on $\partial \Omega$,
and with finite energy semi-norm
\bea
\|u\|_{\mathfrak H_\diamond} = \int_\Omega |\nabla u|^2 dx.
\eea 
Since the functions in  $\mathfrak H_\diamond$ are harmonic
in $D^\prime$, the  
$\mathfrak H_\diamond$
is a closed subspace of $H^1(\Omega)$. Later on, we will give
a new characterization of the space $\mathfrak H_\diamond$ in 
terms of the single layer potential on  $\partial D$ associated with the 
Neumann function of $\Omega$. 
 \\

 We remark that $Tu = 0$ for all
$u$ in $H^1_0(D^\prime)$, and $Tu = u$ for all
$u$ in $H^1_0(D)$ (the set of functions in $H^1(D)$ with trace zero).\\

We also remark that
$T\mathfrak H_\diamond \subset \mathfrak H_\diamond $
and hence the  restriction of $T$ to  $\mathfrak H_\diamond $
defines a linear bounded operator.  Since we are interested in
harmonic functions  in  $D$ and $D^\prime =\Omega\setminus \overline D$,
 we only consider
the action of $T$ on  the closed space $\mathfrak H_\diamond $. We 
further keep the notation $T$ for the restriction of $T$ to 
$\mathfrak H_\diamond $.  We will prove later that $T$ has
only isolated eigenvalues with an accumulation point  $1/2$. 
We denote by $\left(\lambda_n^-\right)_{n\geq 1}$ the eigenvalues of
$T$   repeated according to their multiplicity, and 
 ordered as follows
\bea
 0<\lambda_1^- \leq \lambda_2^- \leq \cdots < \frac{1}{2},
\eea
 in $(0, 1/2]$ and, similarly,
 \bea
1> \lambda_1^+ \geq \lambda_2^+ \geq \cdots > \frac{1}{2}.
\eea
 the eigenvalues in $[1/2, 1)$. The eigenvalue $\lambda_\infty= 1/2$ is the unique accumulation
point of the spectrum.  To ease the notation we further
 denote the orthogonal spectral projector on the eigenspace
associated to  $1/2$, by $\int_{\partial \Omega} \cdot w_\infty^\pm(z) ds(z)w_\infty^\pm(x).$ 
 %%%%%%%%%%%%%%%%%%%%%%%%%%%%%%%%%%%%
Next, we will characterize the spectrum of  $T$ via
the mini-max principle.  
%%%%%%%%%%%%%%%%%%%%%%%%%%%%%%%%%%
\begin{proposition}  The variational Poincar\'e operator has the
following decomposition
\bean \label{Tdecomp}
T= \frac{1}{2}I +K,
\eean
where $K$ is a compact self-adjoint operator.  Let  $w_n^\pm, \; n\geq 1$ 
be the  eigenfunctions associated to the eigenvalues 
$\left(\lambda_n^\pm\right)_{n\geq 0}$. Then 
\bea
\lambda_1^-&=& \min_{0\not= w\in  \mathfrak H_\diamond}
\frac{\int_D|\nabla w(x)|^2 dx}{\int_\Omega |\nabla w(x)|^2 dx},\\
\lambda_n^-&=& \min_{0\not=w \in  \mathfrak H_\diamond, w\perp w_1,\cdots, w_{n-1} }
\frac{\int_D |\nabla w(x)|^2 dx}{\int_\Omega |\nabla w(x)|^2 dx},\\ 
&=& \max_{F_n\subset  \mathfrak H_\diamond, \; dim(F_n) = n-1}\min_{w\in F_n}
\frac{\int_D |\nabla w(x)|^2 dx}{\int_\Omega |\nabla w(x)|^2 dx},\\
\eea

and similarly 
\bea
\lambda_1^+&=& \max_{0\not=w\in  \mathfrak H_\diamond}
\frac{\int_D|\nabla w(x)|^2 dx}{\int_\Omega |\nabla w(x)|^2 dx},\\
\lambda_n^+&=& \max_{0\not=w \in  \mathfrak H_\diamond, w\perp w_1,\cdots, w_{n-1} }
\frac{\int_D |\nabla w(x)|^2 dx}{\int_\Omega |\nabla w(x)|^2 dx},\\ 
&=& \min_{F_n\subset  \mathfrak H_\diamond, \; dim(F_n) = n-1}\max_{w\in F_n}
\frac{\int_D |\nabla w(x)|^2 dx}{\int_\Omega |\nabla w(x)|^2 dx}.\\
\eea
\end{proposition}
%%%%%%%%%%%%%%%%%%%%%%%%%%%%%

%%%%%%%%%%%%%%%%%%%%%%%%%%%%%%%%%%%%%%%

We have the following decomposition of $u(x, \omega)$ in the basis
of the eigenfunctions of the variational Poincar\'e operator $T$.  
%%%%%%%%%%%%%%%%%%%%%%%%%%%%%%%%%%%%%%
%%%%%%%%%%%%%%%%%%%%
\begin{theorem}\cite{AT} \label{freqdepend}
Let $u(x, \omega)$ be the unique solution to the system 
\eqref{maineq}. \\

Then, the following decomposition holds:
\bean \label{decomposition}
u(x, \omega) = k_0^{-1} u_0(x) +   \sum_{n=1}^\infty 
\frac{  \int_{\partial \Omega} f(z) w_n^\pm(z) ds(z) 
}{k_0+\lambda_n^\pm(k(\omega) -k_0) }  w_n^\pm(x),\quad x\in \Omega,
\eean
where $u_0(x) \in H^1_\diamond (\Omega)$ 
depends only on $f$ and $D$, and is the unique solution to 
\bean \label{nonfrequencypart}
\left \{
\ba{lllcc}
\Delta v = 0 & \textrm{in}
\quad D^\prime,\\
\nabla v = 0 & \textrm{in}
\quad D,\\
 \partial_{\nu_\Omega} v = f & \textrm{on}
\quad \partial \Omega.
\ea
\right.
\eean

\end{theorem}
%%%%%%%%%
\proof
In order to have a self-contained  document  we give  the proof of the theorem.\\

 We first observe that frequency dependent part 
$$u_f = u-k_0^{-1} u_0,$$ 
lies in $\mathfrak H_\diamond $. Since the eigenfunctions 
$w^\pm(x)$ form an orthonormal basis of $\mathfrak H_\diamond $, the frequency part
$u_f$ posses the following spectral decomposition:  

\bea
u_f(x) =   \sum_{n=1}^\infty 
 \int_{ \Omega} \nabla u_f(z) \nabla  w_n^\pm(z) dzw_n^\pm(x),  \quad x\in \Omega.
\eea
A forward computation leads to 
\bea
\int_{ \Omega} \nabla u_f(z) \nabla  w_n^\pm(z) dz &= &\int_{ \Omega} \nabla u(z) \nabla  w_n^\pm(z) dz. 
\eea

On the other hand, since $u \in H_\diamond^1(\Omega) $, we obtain
\bea
 \int_{\Omega} \nabla u(z) \nabla  w_n^\pm(z) dz  &=& \lambda_n^\pm  \int_{D} \nabla u(z) \nabla  w_n^\pm(z) dz\\
 &=&\frac{k_0}{k(\omega)}  \lambda_n^\pm  \int_{\partial D} \partial_{\nu_D}
 u(z)|^+ w_n^\pm(z) ds(z)\\
 &=& \frac{k_0}{k(\omega)}  \lambda_n^\pm  
 \int_{D^\prime} \nabla u(z) \nabla  w_n^\pm(z) dz - \frac{k_0}{k(\omega)}  \lambda_n^\pm  
 \int_{\partial \Omega} f(z)w_n^\pm(z) ds(z).
\eea
Using the simple fact that 
\bea
\int_{\Omega} \nabla u(z) \nabla  w_n^\pm(z) dz  &=& \int_{D} \nabla u(z) \nabla  w_n^\pm(z) dz+  \int_{D^\prime} \nabla u(z) \nabla  w_n^\pm(z) dz, 
\eea
we obtain the desired decomposition.

\endproof

%%%%%%%%%%%%%%%%%%%
 In \cite{AT}, assuming that the profile $k(\omega)$ is given, the  spectral decomposition  \eqref{decomposition}  of the solution of the forward conductivity problem has been used to retrieve the Cauchy data corresponding to the extreme case of perfect conductor $u_0$.   Based on unique
 continuation techniques, the uniqueness of the mfEIT problem,  and  rigorous stability estimates 
 have been obtained  from the knowledge of $u_0|_{\partial \Omega}$, in the case where the anomaly is within a class of star shaped domains. \\ 

Assume that $X_0 \in \Omega$, and let   
 $d_1= \textrm{dist}(X_0, \partial \Omega)$ and let $d_0< d_1$.
 For $\delta>0$ small enough, and $m>0$ large enough,  define 
 the set of anomalies:
\bea
\mathfrak D:=  \left\{D =\left\{ X_0+ 
\Upsilon(\theta) 
\begin{pmatrix}
\cos \theta \\ 
\sin \theta
\end{pmatrix},  \theta \in [0;2\pi)\right\};  \Upsilon \in \Pi \right\},
\eea
where 
\bea
\Pi := \left\{ d_0<\Upsilon(\theta)< d_1-\delta;\;\Upsilon(2\pi) = 
\Upsilon(0);\; \|\Upsilon\|_{C^{\beta}} \leq m, \; \beta \geq 2 \right\}.
\eea 

In this paper we consider the  reconstruction of the profile function $k(\omega)$ defined
in \eqref{conductivitydistribution}, and anomalies $D$  within the set $\mathfrak D$.   
At  first glance,  the numerical reconstruction  of the inclusion in the mfEIT problem does not  
need to follow the path of the  theoretical results derived in  \cite{AT}, that is, to determine first 
$u_0|_{\partial \Omega}$, and then find the high conductor anomaly $D$ that produces 
the recovered  Cauchy data  of $u_0$.  In fact preliminary numerical calculations show
that  a blind minimization  approach that searches the anomaly  $D$ and the  profile 
$k(\omega)$  using boundary multifrequency data does not converge in most cases, and
if it happens to converge the rate turns out to be very slow. These difficulties  are well known in
inverse  conductivity problem, usually it is very hard to distinguish between the
 conductivity value and the size of the  anomaly \cite{AK}.  On the other hand 
 the numerical identification of  a high conductor anomaly  is  a well known inverse problem,
 and many works have been done on it  (see for instance \cite{KS, LL}).   We can cite for example
quasi-reversibility type based methods  \cite{KS, LL}.  Here, we will consider the parameterization 
type based methods  \cite{CK,ACLZ,Ru}. Since the problem is ill-posed we will use  a cut-off
regularization approach that  consists on  taking into account in the computation only the important 
Fourier modes of the parametrization function 
$\Upsilon(\theta)$. Then, the identification is transformed into an optimization problem with
a finite number of degree of freedom.  The problem is still strongly nonlinear
we propose here to solve it using the gradient  adjoint 
method.  \\

The algorithm we propose in this paper for identifying numerically the anomaly
and the frequency conductivity profile  is inspired by the theoretical approach
developed in\cite{AT}, and it can be summarized as follows:\\

(i)To recover $u_0(x)|_{\partial \Omega}$ and $\kappa_p, \; p=1, 2, 3,$  from the knowledge 
of $u(x, \omega)|_{\partial \Omega},\; \omega \in (\underline \omega, \overline \omega)$.  Here,
we will use an linear algebraic approach based on the understanding of
 the behavior of the spectrum of Neumann-Poincar\'e operator near its unique accumulation point, 
 and the sparsity of the considered conductivity distribution. \\
 
(ii) To identify the anomaly $D$ from the Cauchy data $(u_0(x), f(x))$
on the boundary $\partial \Omega$ using a cut-off parameterization/Fourier approach. 
Here  to further stabilize and speed up the 
convergence of   the iterative gradient based method we will use two linearly independent  boundary currents $(f_1, f_2)$. 

%%%%%%%%%%%%%%%%%%%%%%%%%%%% 
%%%%%%%%%%%%%%%%%%%%%%%%%%%%
\section{The linearized map}
%%%%%%%%%%%%%%%%%%%%%%%%%%%%%
We shall use gradient methods to  identify the anomaly from  the Cauchy 
data $(u_0(x), f(x))$ on the boundary $\partial \Omega$. In this section, we determine 
$u_h$,  the derivative of $u_0$ with respect to a shape perturbation in the direction
 $h(x)\nu_D(x)$, where $h(x) $ is a scalar function defined on $\partial D$. We will follow
 the analysis in  \cite{AKLZ} for a non-degenerate conductivity inside the anomaly, based on integral
 equations techniques.  \\
 
Let $D\in \mathfrak D$ be a given anomaly. We define
 $X(t):[a,b]\rightarrow \mathbb{R}^2$ to be a smooth 
 clockwise parametrization of $\partial D$, where 
 $a,b \in \mathbb{R}$, $a<b$. We assume that 
 $X \in \mathcal{C}^{\beta}([a,b])$ and $|X^\prime(t)|=1$ for all $t\in[a,b]$. Then
 \begin{equation}
 \partial D=\{x= X(t), t\in [a,b]\}.
\end{equation}
Let $h\in C^2(\partial D)$, and define 
 the boundary of the perturbed anomaly $D_\varepsilon$  by
\begin{equation}
 \partial D_\varepsilon=\{\tilde{x}=\tilde{X}(t):=X(t)+\varepsilon h(X(t)) \nu_D(X(t)), t\in [a,b]\}.
\end{equation}

Define $u_\varepsilon$  to be the unique solution to the system 
\eqref{nonfrequencypart} associated
to the perturbed anomaly  $D_\varepsilon$, that is 
   \begin{equation}\label{Pbper}
   \left\{
    \begin{array}{lr}
  \bigtriangleup v =0 & \text{in } \Omega \setminus \overline{D_\varepsilon}, \\
  \nabla v = 0 & \text{in } D_\varepsilon,\\
\partial_{\widetilde \nu } v =f & \text{on } \partial \Omega,\\
  \int_{\partial \Omega} v d\sigma =0,
   \end{array}
  \right.
   \end{equation}
   where $\widetilde \nu $ is the outward normal vector on $\partial D_\varepsilon$.\\

The objective of this section is to derive a linear correction $u_h$
of $u_\varepsilon$,   such that 
\bea
u_\varepsilon =  u_0+ \varepsilon u_h +O(\varepsilon^2), \qquad \textrm{ as } \varepsilon \to 0. 
\eea

The main result of this section is the following.  
%%%%%%%%%%%%%%%%%%%%%%%%%%%%
\begin{theorem} \label{frechet} Let  $h\in C^2(\partial D)$ be fixed. Then, $u_h$
is the unique solution to the system 
 \begin{equation}\label{Pbuh}
   \left\{
    \begin{array}{lr}
  \bigtriangleup v =0 & \text{in } \Omega \setminus \overline{D}, \\
  \nabla v = 0 & \text{in } D,\\
u|_+-u|_-=-h \partial_{\nu_D} u_0|_+ & \text{on } \partial D,\\
\partial_{\widetilde \nu } v = 0 & \text{on } \partial \Omega,\\
 \int_{\partial \Omega} v d\sigma =0.
   \end{array}
  \right.
   \end{equation}

\end{theorem}
%%%%%%%%%%%%%%%%%%%%%%%

\proof
We first derive an integral equation representation  of the field.  \\

Let $G(x,z)=  \frac{1}{2\pi}\log(|x-y|)$ 
be the Green function for the Laplacian  in
$\mathbb R^2$, and define the single layer potentials respectively 
on $\partial D$, and $\partial \Omega$
 by   \\
\bea
 \mathcal S_D: H^{-\frac{1}{2}}(\partial D)\rightarrow H^{\frac{1}{2}}(\partial D),\\
  \mathcal S_D \varphi(x) = \int_{\partial D} G (x, z)  \varphi(z) ds(z),
\eea
and
\bea
 \mathcal S_\Omega: H^{-\frac{1}{2}}(\partial \Omega)\rightarrow 
H^{\frac{1}{2}}(\partial \Omega),\\
  \mathcal S_\Omega \psi(x) = \int_{\partial \Omega} G (x, z)  \psi(z) ds(z),
\eea

 They  satisfy the following jump relations through 
 the boundary of  respectively $D$ and $\Omega$  \cite{AK}
 \bea
 \partial_{\nu_D} S_D\varphi(x)|^\pm = \pm \frac{1}{2}\varphi(x) + \mathcal K^*_D\varphi(x), 
 \qquad \textrm{for }  x\in \partial D,\\
 \partial_{\nu_\Omega} S_\Omega \varphi(x)|^\pm = \pm \frac{1}{2}\varphi(x) 
 + \mathcal K^*_\Omega\varphi (x),  \qquad \textrm{for }  x\in \partial \Omega,
 \eea
  where 
 \bea
 \mathcal K^*_D:
   H^{-\frac{1}{2}}(\partial D)\rightarrow H^{-\frac{1}{2}}(\partial D)\\
 \mathcal K^*_D\varphi(x) = \int_{\partial D} \partial_{\nu_D(x)}  
G(x, z)  \varphi(z) ds(z), 
 \eea
 and
 \bea
 \mathcal K^*_\Omega:
   H^{-\frac{1}{2}}(\partial \Omega)\rightarrow H^{-\frac{1}{2}}(\partial \Omega)\\
 \mathcal K^*_\Omega\psi (x) = \int_{\partial \Omega} \partial_{\nu_\Omega(x)}  
G(x, z)  \psi(z) ds(z), 
 \eea
 are  compact operators.  \\
 
 Since $G(x,z)$ is harmonic in $\mathbb R^2\setminus \{z\}$, 
 $S_\Omega\psi$ has a unique  harmonic extension in
 $\mathbb R^2 \setminus \partial \Omega$.  Similarly, $S_D \varphi$ has a 
 unique  harmonic extension in
 $\mathbb R^2 \setminus \partial D$.  In addition, we have 
 \bea
 \mathcal S_D \varphi |^+=   \mathcal S_D \varphi |^- \qquad \textrm{on } \partial D,\\
  \mathcal S_\Omega \psi |^+=   \mathcal S_\Omega \psi |- \qquad \textrm{on } \partial \Omega.
 \eea
We also define the double layer potential $\mathcal D_D $ 
 \bea
 \mathcal D_D: H^{-\frac{1}{2}}(\partial D)\rightarrow H^{\frac{1}{2}}_{loc}(\mathbb R^2\setminus
 \partial D),\\
  \mathcal D_D \varphi(x) = \int_{\partial D} \partial_DG (x, z)  \varphi(z) ds(z).
\eea
It satisfies the following jump relations 
\bea
\partial_{\nu_D}\mathcal D_D \varphi(x) |^+=\partial_{\nu_D}\mathcal D_D \varphi(x) |^-\qquad \textrm{for }  x\in \partial D,\\
 \mathcal D_D \varphi (x)|^\pm =  \pm \frac{1}{2}\varphi(x) + \mathcal K_D\varphi(x), 
 \qquad \textrm{for }  x\in \partial D,\\
\eea
where 
\bea
 \mathcal K_D:
   H^{-\frac{1}{2}}(\partial D)\rightarrow H^{-\frac{1}{2}}(\partial D)\\
 \mathcal K_D\varphi(x) = \int_{\partial D} \partial_{\nu_D(z)}  
G(x, z)  \varphi(z) ds(z), 
 \eea
is the $L^2(\partial D)$-adjoint of $\mathcal K_D^*$.\\

The solution $u_0$ can be written as
  \begin{equation}\label{representation}
   u_0(x)=\mathcal{S}_D \phi (x)+\mathcal{S}_\Omega \psi(x) \qquad x\in \Omega,
  \end{equation}
  where $\phi= \partial_{\nu_D}u_0 |_+  \in H^{-\frac{1}{2} }_\diamond (\partial D)$ 
  and $\psi \in H^{-\frac{1}{2}}(\partial \Omega)$.\\

  Similarly, we have 
  \begin{equation}
   u_\varepsilon(x)=\mathcal{S}_{D_\varepsilon} \phi_\varepsilon (x)+\mathcal{S}_\Omega \psi_\varepsilon(x),  \qquad x\in \Omega,
  \end{equation}
  where $\phi_\varepsilon= \partial_{\tilde \nu}u_\varepsilon |_+  \in H^{-\frac{1}{2} }_\diamond (\partial D_\varepsilon)$ 
  and $\psi_\varepsilon \in H^{-\frac{1}{2}}(\partial \Omega)$.\\

  Using the jump relations and the facts that $  \partial_{\nu_D}u_0 |_- =0$ 
  on $\partial D$, and $\partial_{\nu_\Omega } u_0 |_-=f$ on $\partial \Omega$, the densities 
  $\phi$ and $\psi$ satisfy the following system
  \bean \label{Sys0}
  (-\frac{1}{2} I +\mathcal{K}^*_D)\phi(x)+ \partial_{\nu_D}\mathcal{S}_\Omega \psi(x) &=0 & \text{on } \partial D, \\
 \partial_{\nu_\Omega} \mathcal{S}_D \phi(x) +(-\frac{1}{2}I+\mathcal{K}^*_\Omega)\psi(x) &=f & \text{on } \partial \Omega.
\eean
 This system can be also represented in a matrix form
 \bean
  \label{matrix_origin}
 M
 \begin{pmatrix}
  \phi \\
  \psi
 \end{pmatrix}
:=
  \begin{pmatrix}
  -\frac{1}{2} I +\mathcal{K}^*_D &  \partial_{\nu_D}\mathcal{S}_\Omega \\
  \partial_{\nu_\Omega }\mathcal{S}_D & -\frac{1}{2}I+\mathcal{K}^*_\Omega
  \end{pmatrix}
  \begin{pmatrix}
   \phi \\
   \psi
  \end{pmatrix}
  =
  \begin{pmatrix}
   0 \\
   f
  \end{pmatrix}.
 \eean
 
The same analysis leads to the system 
\begin{equation}\label{matrix_pertub}
 M_\varepsilon
 \begin{pmatrix}
  \phi_\varepsilon \\
  \psi_\varepsilon
 \end{pmatrix}
:=
  \begin{pmatrix}
  -\frac{1}{2} I +\mathcal{K}^*_{D_\varepsilon} &  
  \partial_{\widetilde \nu}\mathcal{S}_\Omega \\
  \partial_{\nu_\Omega} \mathcal{S}_{D_\varepsilon} & -\frac{1}{2}I+\mathcal{K}^*_\Omega
  \end{pmatrix}
  \begin{pmatrix}
   \phi_\varepsilon \\
   \psi_\varepsilon
  \end{pmatrix}
  =
  \begin{pmatrix}
   0 \\
   f
  \end{pmatrix}.
\end{equation}

From the parameterization of $\partial D$, we deduce that 
the outward unit normal vector $\nu_D(x)$ is given by $\nu_D(X(t))=R_{-\frac{\pi}{2}}T(X(t))$, where $R_{-\frac{\pi}{2}}$ is the rotation with the angle $-\frac{\pi}{2}$, and  
 $T(X(t))=X^\prime(t)$ is the tangential normal vector.  Let  $\gamma(X(t))$ be the curvature, it 
 satisfies
\begin{equation}
X^{\prime\prime}(t)=\gamma(X(t))\nu(X(t)). 
\end{equation}
Using the parameterization of  $\partial D_\varepsilon$, we deduce the following
asymptotic expansion 
\begin{equation}\label{asym_nu}
 \tilde{\nu}(\tilde{x})=\nu(x)-\varepsilon h^\prime(t)T(x)+O(\varepsilon^2), \qquad \textrm{for } \tilde x\in 
 \partial D_\varepsilon,
\end{equation}
where $h^\prime(t)=\frac{d}{dt}h(X(t))$ (we also use $h^\prime (x)$ 
to denote this quantity).
In the same way obtain the asymptotic expansion of the
 length element 
\begin{equation}\label{asym_ds}
 ds_\varepsilon(\tilde{z})=ds(z)(1-\varepsilon\gamma(y)h(z)+O(\varepsilon^2)).
\end{equation}
Let $\Psi_\varepsilon$ be the diffeomorphism  from $\partial D$ onto $\partial D_\varepsilon$ given by 
$\Psi_\varepsilon(x)=x+\varepsilon h(x)\nu(x)$. From \cite{AKLZ}, we deduce the asymptotic expansion of $\mathcal{K}^*_{D_\varepsilon}$
\begin{equation}
 (\mathcal{K}^*_{D_\varepsilon}\tilde{\phi})\circ \Psi_\varepsilon = \mathcal{K}^*_D \phi + \varepsilon\mathcal{K}^{(1)}_D \phi + O(\varepsilon^2),
\end{equation}
where $\tilde{\phi}=\phi \circ \Psi^{-1}_\varepsilon$,  and the operator $\mathcal{K}^{(1)}_D$ is defined by
\begin{align}
 \mathcal{K}^{(1)}_D \phi (x)= \frac{1}{2\pi}\int_{\partial D} [(&\frac{1}{|x-y|^2}-\frac{2\langle x-y, \nu_D(x) \rangle^2}{|x-y|^4})h(x)-\frac{\langle x-y,T(x) \rangle}{|x-y|^2}h'(x) \nonumber \\
 &-\frac{\langle\nu_D(x),\nu_D(y)\rangle}{|x-y|^2}h(y)+\frac{2\langle x-y,\nu_D(x) \rangle \langle x-y, \nu_D(y) \rangle}{|x-y|^4}h(y) \nonumber \\
 &-\frac{\langle x-y, \nu_D(x)\rangle}{|x-y|^2}\gamma(y)h(y)]\phi(y)ds(y)
\end{align}
Now, we calculate the asymptotic expansion of the operators $\partial_{\tilde \nu}\mathcal{S}_\Omega$ on $\partial D_\varepsilon$ and $\partial _{\nu_\Omega}\mathcal{S}_D$ on $\partial \Omega$.\\

Let $\psi \in H^{-\frac{1}{2}}(\partial \Omega)$ be fixed.   Using the relation 
$\tilde{x}=x+\varepsilon h(x) \nu(x) \in \partial D_\varepsilon$ for  $x \in \partial D$,
we obtain 
\begin{align*}
 \partial_{\tilde \nu}  \mathcal{S}_\Omega \psi (\tilde{x})&= \frac{1}{2\pi}\int_{\partial \Omega}\frac{\langle \tilde{x}-y,\tilde{\nu}(\tilde{x}) \rangle}{|\tilde{x}-y|^2}\psi(y)ds(y),\\
 &=\frac{1}{2\pi}\int_{\partial \Omega}\frac{\langle x+\varepsilon h(x) \nu_D(x)-y,\nu_D(x)-\varepsilon h^\prime(x)T(x) \rangle}{|x+\varepsilon h(x) \nu(x)-y|^2}\psi(y)ds(y)+O(\varepsilon^2),\\
 &= \partial_{\nu_D}\mathcal{S}_\Omega \psi (x)+\varepsilon(-h^\prime(x)
 \partial_T\mathcal{S}_\Omega \psi (x)+h(x) \mathcal{S}^{(1)}_\Omega \psi(x))+O(\varepsilon^2),
\end{align*}
where $\partial_T$ denotes the tangential derivative, and $\mathcal{S}^{(1)}_\Omega$ is defined by
\begin{equation}
 \mathcal{S}^{(1)}_\Omega \psi(x)= \frac{1}{2\pi}\int_{\partial \Omega} [\frac{1}{|x-y|^2}-\frac{2\langle x-y,\nu(x) \rangle^2}{|x-y|^4}]\psi(y) ds(y),
\end{equation}
for $x \in \partial D$.\\

We  further determine the asymptotic expansion of $\partial_{\nu_\Omega} 
\mathcal{S}_{D_\varepsilon}$ on $\partial \Omega$.
Let $\phi \in H^{-\frac{1}{2}}(\partial D)$,  and $x \in \partial \Omega$, we have
\begin{align}
 &\partial_{\nu_\Omega}\mathcal{S}_{D_\varepsilon} \tilde{\phi} (x)= \frac{1}{2\pi}\int_{D_\varepsilon}\frac{\langle x-\tilde{y}, \nu_\Omega(x) \rangle}{|x-\tilde{y}|^2}\tilde{\phi}(\tilde{y})ds_\varepsilon (\tilde{y}),\nonumber\\
 &= \frac{1}{2\pi}\int_D\frac{\langle x-y-\varepsilon h(y) \nu_D(y), \nu_\Omega(x) \rangle}{|x-y-\varepsilon h(y) \nu_D(y)|^2}\phi(y))(1-\varepsilon\gamma(y)h(y))ds+ O(\varepsilon^2),\nonumber\\
 &=\partial_{ \nu_\Omega} \mathcal{S}_D \phi (x)+\varepsilon \{ \frac{1}{2\pi}
 \int_{\partial D}[-\frac{\langle \nu_D(y),\nu_\Omega(x) \rangle}{|x-y|^2}+2\frac{\langle x-y, \nu_\Omega(x) \rangle \langle x-y, \nu_D(y) \rangle}{|x-y|^4}]h(y) \phi(y)ds(y)\nonumber\\
 &-\frac{1}{2\pi}\int_{\partial D}\frac{\langle x-y, \nu_\Omega(x) \rangle}{|x-y|^2}\gamma(y)h(y)\phi(y)ds(y)\}+O(\varepsilon^2),\nonumber \\
 &= \partial_{ \nu_\Omega} \mathcal{S}_D \phi (x)+\varepsilon \partial_{\nu_\Omega}\left(\mathcal{D}_D (h\phi)-\mathcal{S}_D(\gamma h \phi)\right)(x)+O(\varepsilon^2).
\end{align}
Consequently
\bea
M_\varepsilon=M+\varepsilon M_h + O(\varepsilon^2),
\eea 
where the operator $M_h$ on $H^{-\frac{1}{2}}_\diamond(\partial D) \times 
H^{-\frac{1}{2}}(\partial \Omega)$ is defined by
\begin{equation}
 M_h :=
  \begin{pmatrix}
   \mathcal{K}^{(1)}_D & -h^\prime\partial_T\mathcal{S}_{\Omega}+h\mathcal{S}^{(1)}_\Omega\\
   \partial_{\nu_\Omega}[\mathcal{D}_D(h\cdot)-\mathcal{S}_D(\gamma h \cdot)] & 0
  \end{pmatrix}.
\end{equation}

So, the systems (\ref{matrix_origin}) and (\ref{matrix_pertub}) imply 
\begin{equation}
 \begin{pmatrix}
  \phi_\varepsilon \\
  \psi_\varepsilon
 \end{pmatrix}
=
\begin{pmatrix}
  \phi \\
  \psi
 \end{pmatrix}
 +
 \varepsilon \begin{pmatrix}
  \phi_h\\
  \psi_h
 \end{pmatrix}+O(\varepsilon^2),
\end{equation}
where $\begin{pmatrix}
  \phi_h\\
  \psi_h
 \end{pmatrix}$ is given by 
 \begin{equation}\label{relation_matrix_phi}
  \begin{pmatrix}
  \phi_h\\
  \psi_h
 \end{pmatrix}=-M^{-1}M_h \begin{pmatrix}
  \phi\\
  \psi
 \end{pmatrix}.
 \end{equation}
Thus, using the representation formula and following the same calculus, we   
determine the asymptotic expansion of the solution $u_\varepsilon|_{\partial \Omega}$
\begin{equation}
 u_\varepsilon(x)=u_0(x)+\varepsilon(\mathcal{S}_D \phi_h(x)+\mathcal{S}_\Omega \psi_h(x)+\mathcal{D}_D(h\phi)(x)-\mathcal{S}_D(\gamma h \phi)(x))+O(\varepsilon^2).
\end{equation}
We further denote
\bea
\tilde{u}_h&=& \mathcal{S}_D \phi_h+\mathcal{S}_\Omega \psi_h,\\ 
u_h&=&\tilde{u}_h+\mathcal{D}_D(h\phi)-\mathcal{S}_D(\gamma h \phi).
\eea

We deduce from  (\ref{relation_matrix_phi})
\begin{equation}\label{egalite1}
 \partial_{\nu_D } \tilde{u}_h|_- + \mathcal{K}^{(1)}_D \phi -h^\prime \partial_T \mathcal{S}_{\Omega}\psi+h\mathcal{S}^{(1)}_\Omega\psi =0,
\end{equation}
on $\partial D$,
and 
\begin{equation}\label{egalite2}
 \partial_{\nu_\Omega} \tilde{u}_h|_- +   \partial_{\nu_\Omega} 
 [\mathcal{D}_D(h\phi)-\mathcal{S}_D(\gamma h \phi)]=0,
\end{equation}
on $\partial \Omega$, and hence 
\begin{equation}
\partial_{\nu_\Omega} u_h|_-=0 \qquad \text{on } \partial \Omega.
\end{equation}
The equality (\ref{representation}) and the fact that $\partial_T u_0 =0$ on $\partial D$, lead to
\begin{align*}
 0&= \partial_T\left( h \partial_T u_0(x)\right),\\
 &=\partial_T \left(h \partial_T \mathcal{S}_D\phi(x)\right)+h^\prime(x)\partial_T\mathcal{S}_\Omega \psi(x)+h(x)\partial_T^2\mathcal{S}_\Omega \psi(x),\\
 &=\partial_T \left(h \partial_T \mathcal{S}_D\phi(x)\right)+h^\prime(x)\partial_T\mathcal{S}_\Omega \psi(x)\\
 &+h(x)\frac{1}{2\pi}\int_{\partial \Omega}[\frac{-1}{|x-y|^2}+2\frac{\langle x-y,\nu_D(x)\rangle^2}{|x-y|^4}+\gamma(x)\frac{\langle x-y,\nu_D(x) \rangle}{|x-y|^2}]\psi(y)ds(y),\\
 &= \partial_T (h \partial_T \mathcal{S}_D\phi(x))+h^\prime(x)\partial_T\mathcal{S}_\Omega \psi(x)-h(x)\mathcal{S}^{(1)}_\Omega \psi(x)+\gamma(x)h(x) \partial_{\nu_D}\mathcal{S}_\Omega \psi(x),
\end{align*}
which implies
\begin{align*}
 -h^\prime \partial_T\mathcal{S}_{\Omega}\psi+h\mathcal{S}^{(1)}_\Omega\psi=\partial_T(h\partial_T\mathcal{S}_D\phi)+\gamma h(\frac{1}{2}I-\mathcal{K}^*_D)\phi.
\end{align*}
A similar calculus gives
\begin{align*}
 &\partial_T (h \partial_T \mathcal{S}_D \phi (x)) \nonumber\\
 &= \frac{1}{2\pi}\int_{\partial D}[h'(x)\frac{\langle x-y, T(x) \rangle}{|x-y|^2}+h(x)(\frac{-1}{|x-y|^2}+\frac{2\langle x-y, \nu_D(x) \rangle^2}{|x-y|^4}+2\gamma(x)\frac{\langle x-y, \nu_D(x) \rangle}{|x-y|^2})]\phi(y)ds(y),
\end{align*}
for $x \in \partial D$.\\

Thus
\begin{align}\label{simplify_KD1}
 &\mathcal{K}^{(1)}_D\phi(x)-h^\prime \partial_T\mathcal{S}_{\Omega}\psi(x)+h\mathcal{S}^{(1)}_\Omega\psi(x) \nonumber\\
 &= \frac{1}{2\pi}\int_{\partial D}[-\frac{\langle \nu_D(x),\nu_D(y) \rangle}{|x-y|^2}+\frac{\langle x-y, \nu_D(x) \rangle \langle x-y, \nu_D(y) \rangle}{|x-y|^4}]h(y) \phi(y)ds(y)\nonumber\\
 &-\frac{1}{2\pi}\int_{\partial D}\frac{\langle x-y, \nu_D(x) \rangle}{|x-y|^2}\gamma(y)h(y)\phi(y)ds(y)+\frac{1}{2}\gamma(x)h(x)\phi(x).
\end{align}
By the continuity of the normal derivative of double layer potentials and the jump relation, we have, for $x\in \partial D$,
\begin{align}\label{DDSD_dD}
& \partial_{\nu_D}[\mathcal{D}_D(h\phi)-\mathcal{S}_D(\gamma h \phi)](x)|_{-}\nonumber\\
& =\frac{1}{2\pi}\int_{\partial D}[-\frac{\langle \nu_D(x),\nu_D(y) \rangle}{|x-y|^2}+2\frac{\langle x-y, \nu_D(x) \rangle \langle x-y, \nu_D(y) \rangle}{|x-y|^4}]h(y) \phi(y)ds(y)\nonumber\\
&+\frac{1}{2}\gamma(x)h(x)\phi(x)-\frac{1}{2\pi}\int_{\partial D}\frac{\langle x-y, \nu_D(x) \rangle}{|x-y|^2}\gamma(y)h(y)\phi(y)ds(y).
\end{align}

Using (\ref{egalite1}), (\ref{simplify_KD1}) and (\ref{DDSD_dD}), we have,
\bea
 \partial_{\nu_D} u_h|_- = \partial_{\nu_D} \tilde{u}_h|_-+ \partial \nu_D[\mathcal{D}_D(h\phi)-\mathcal{S}_D(\gamma h \phi)]|_{-} =0,
\eea
which gives the desired result. 
\endproof

\section{Reconstruction of $u_0(x)|_{\partial \Omega}$ and  $k(\omega)$}
%%%%%%%%%%%%%%%%%%%%%%%%%%%%%%%%%%%%%%%%%%%%%%
%%%%%%%%%%%%%%%%%%%%%%%%%%%
In this section we construct $u_0|_{\partial \Omega}$ from the knowledge
of $u(x, \omega)|_{\partial \Omega},\; \omega \in (\underline \omega, \overline \omega)$.
Here we recall the spectral decomposition  \eqref{decomposition}, also valid
on the boundary $\partial \Omega$, which is the keystone
of our approach. 
\bea
u(x, \omega) = k_0^{-1} u_0(x) +   \sum_{n=1}^\infty 
\frac{  \int_{\partial \Omega} f(z) w_n^\pm(z) ds(z) 
}{k_0+\lambda_n^\pm(k(\omega) -k_0) }  w_n^\pm(x), \quad x\in \partial \Omega.
\eea

The first observation is that the functions $w_n^\pm(x)$ do not need to be orthogonal 
on the boundary $\partial \Omega$. Then, varying the frequency, and so the 
coefficients of the expansion above do not guarantee  the complete  separation 
between the frequency  and the non frequency parts.  The second observation is that  the 
simultaneous  determination of the plasmonic resonances $\lambda_n^\pm, \,n \geq
1$, the frequency profile $k(\omega)$, and $u_0|_{\partial \Omega}$ is strongly nonlinear
while if we assume that $k(\omega)$ and  $\lambda_n^\pm, \,n \geq 1$ are given the problem 
becomes a linear one. \\

We further consider $M\geq 2$  frequencies  $\omega_1, \cdots, \omega_M$ in
$ (\underline \omega, \overline \omega)$,
 and their
associated solutions $u(x,\omega_1), \cdots, u(x,\omega_M)$. 
Since $1/2$ is the unique accumulation point of the eigenvalues 
$(\lambda_n^\pm)_{n\geq 1}$, we only consider the $N_f\geq 0$
 first eigenvalues as unknown variables, and we approximate the others
  eigenvalues by the limiting value  $1/2$. In fact it has been shown in \cite{MS}
  that if $D$ is $C^\beta$  with $\beta\geq 2$ then for  any $\alpha >-2\beta +3,$
  we have
  \bea
 | \lambda_n^\pm -1/2| =  o(n^{\alpha}), \qquad n \to +\infty.
  \eea
  
Thus the boundary regularity is essential to the decay rate of eigenvalues. 
 Consequently if the boundary is $C^\infty$ smooth, then the plasmonic 
 eigenvalues  will decay
 faster than any power order.   Recently  H. Kang and his collaborators have proved 
 the exponential convergence of  the  eigenvalues  in the case of analytic anomalies  \cite{AKM}. 
This theoretical work justifies the exponential decay behavior  that has been observed  
numerically \cite{PP}, and checked   
for sample  geometries  like ellipses. If  $\alpha>0$ is   the modified maximal Grauert radius of
$\partial D$, then
\bea
 | \lambda_n^\pm -1/2| =  O(e^{-n\alpha}), \qquad n \to +\infty.
  \eea
These asymptotic  properties of the spectrum of the Neumann-Poincar\'e operator
suggest to consider only a finite number of them in the spectral decomposition  
 \eqref{decomposition}. We further make the following approximation  
   for $x\in \Omega$, $1\leq p \leq M:$
\bean \label{decom_appro}
u(x,\omega_p)\approx \frac{2}{k(\omega_p)+k_0}v_{1}(x), \qquad \textrm{if  } N_f=0,\\
u(x,\omega_p)\approx k_0^{-1}u_0(x)+\sum_{n=1}^{N_f} \frac{1}{k_0+
\lambda^\pm_n(k(\omega_p)-k_0)}v^\pm_n(x)+\frac{2}{k(\omega_p)+k_0}v_{N_f+1}(x), \qquad \textrm{if  } N_f\geq 1, \nonumber
\eean
where \bea
v^\pm_n (x)=\int_{\partial \Omega} f(z)w^\pm_n (z) ds(z)w^\pm_n (x), \;
v_{N_f+1}(x)=\sum_{n > N_f}\int_{\partial \Omega} f(z)w^\pm_n (z) ds(z)w^\pm_n (x).
\eea
A simple integration by parts, leads for all $n \geq 1$
\begin{equation}
\int_{\partial \Omega} f(z)w^\pm_n (z)ds(z)=\int_\Omega \nabla \mathfrak{f}(x)\nabla w^\pm_n(x)dx,
\end{equation}
where $\mathfrak{f}$ is the unique solution in $H^1_\diamond(\Omega)$ to
\begin{equation}
\left\{
    \begin{array}{lr}
  \bigtriangleup \mathfrak{f}=0  & \text{in } \Omega, \\
  \partial_\nu \mathfrak{f}=f & \text{on } \partial \Omega.
   \end{array}
  \right.
\end{equation}
Consequently, the function 
\bea
 \mathcal P_0 \mathfrak{f}&=&\sum_{n=1}^{N_f+1} v^\pm_n(x),
 \eea
where $ \mathcal P_0$ is the orthogonal projection  onto  the space $\mathfrak{H}_\diamond$.  On the other hand,  $u_0$ satisfies
\begin{align}
&\int_\Omega \nabla u_0(x) \nabla w^\pm_n(x) dx=\int_{\Omega \setminus \overline{D}} \nabla u_0(x) \nabla w^\pm_n(x) dx \nonumber \\
&=\int_{\partial \Omega} u_0(x)
\partial_ {\nu_\Omega}w^\pm_n (x) ds(x)-
\int_{\partial D} u_0(x) \partial_ {\nu_D} w^\pm_n (x) ds(x) = 0,\nonumber
\end{align}
for all $n \geq 1$. \\

Since $ \mathcal P_0 u_0 = 0$, and $ \mathfrak{f}-u_0 \in \mathfrak{H}_\diamond$  the orthogonal projection of $\mathfrak{f}$ onto the space $\mathfrak{H}_\diamond$ is $\mathfrak{f}-u_0$, that is $ \mathfrak{f}-u_0= \mathcal P_0 (\mathfrak{f}-u_0)   =  
 \mathcal P_0 \mathfrak{f}.$ \\

Therefore,  the formula (\ref{decom_appro}) becomes
\bean \label{1mode}
u(x,\omega_p) \approx   \frac{k(\omega_p)-k_0}{k_0(k(\omega_p)+k_0)}u_0(x)+
\frac{2}{k(\omega_p)+k_0}\mathfrak{f}(x), \nonumber \\
 \qquad \textrm{if  } N_f=0, \nonumber\\
\approx \frac{k(\omega_p)-k_0}{k_0(k(\omega_p)+k_0)}u_0(x)+\frac{2}{k(\omega_p)+k_0}\mathfrak{f}(x) 
+\sum_{n=1}^{N_f}(\frac{1}{k_0+\lambda^\pm_n(k(\omega_p)-k_0)}-\frac{2}{k(\omega_p)+k_0})v^\pm_n(x),
\nonumber
\\ \qquad \textrm{if  } N_f\geq 1.
\eean
%%%%%%%%%%%%%%--------------------
Next, we reconstruct $\kappa_1$, $\kappa_2$, $\kappa_3$ and $u_0(x)$ by an optimization algorithm. In order to do so, we need an apriori estimation of the eigenvalues $\widetilde{\lambda_n^\pm}\in [0,1]$ for $ n= 1, \cdots, N_f$. Since the eigenvalues $\lambda_n^\pm$ are within a relative narrow interval, preliminary 
calculations showed that  the reconstruction of $u_0$ is indeed not very  sensitive to the choice of those eigenvalues.  In the rest of this section, we assume that we have we fix $\widetilde{\lambda_n^\pm}\in [0,1]$ for $ n= 1, \cdots, N_f$.\\

Let $(x_j)_{1\leq j \leq N_d} \in \partial \Omega$ a discretization of the boundary $\partial \Omega$, and
 define, for $ n= 1, \cdots, N_f$ the scalar functionals
\begin{align}\label{def_Fj}
&F_j(U_0^{(j)},V^{\pm (j)}_1,\cdots,V^{\pm (j)}_{N_f},\omega,\kappa_1,\kappa_2,\kappa_3):=\nonumber \\
&\frac{k(\omega,\kappa_1,\kappa_2,\kappa_3)-k_0}{k_0(k(\omega,\kappa_1,\kappa_2,\kappa_3)+k_0)}U_0^{(j)}+\frac{2}{k(\omega,\kappa_1,\kappa_2,\kappa_3)+k_0}\mathfrak{f}(x_j)\nonumber \\ 
&+\sum_{n=1}^{N_f}(\frac{1}{k_0+\widetilde{\lambda^\pm_n}(k(\omega,\kappa_1,\kappa_2,\kappa_3)-k_0)}-\frac{2}{k(\omega,\kappa_1,\kappa_2,\kappa_3)+k_0})V^{\pm (j)}_n.
\end{align}
where $(U_0^{(j)})_{1\leq j \leq N_d}$ and $(V_n^{\pm(j)})_{1\leq j \leq N_d}$ are
 vectors in $\mathbb{R}^{N_d}$, that approximate respectively $(u_0(x_j))_{1\leq j \leq N_d}$ and $(v_n(x_j))_{1\leq j \leq N_d}$. \\

The scheme consists in  minimizing the scalar functional
\begin{align}
&J_m(U_0,V^\pm_1,\cdots,V^\pm_{N_f},\kappa_1,\kappa_2,\kappa_3):=\nonumber \\
&\frac{1}{2}\sum_{p=1}^M \sum_{j=1}^{N_d} |u(x_j,\omega_p)-F_j(U_0^{(j)},V^{\pm (j)}_1,\cdots,V^{\pm (j)}_{N_f},\omega_p,\kappa_1,\kappa_2,\kappa_3)|^2.
\end{align}
So, we can easily calculate its gradient from (\ref{def_Fj}) and (\ref{empirical}), for $i=1,2,3$, $1\leq l \leq N_d$ and $1\leq n\leq N_f$,
\begin{align}
\frac{\partial J_m}{\partial \kappa_i} = \sum_{p=1}^M \sum_{j=1}^{N_d}(\overline{u(x_j,\omega_p)-F_j(\cdot,\omega_p)})\frac{\partial F_j}{\partial \kappa_i}(\cdot,\omega_p), \label{dJmdki}\\
\frac{\partial J_m}{\partial U_0^{(l)}}= \sum_{p=1}^M (\overline{u(x_j,\omega_p)-F_l(\cdot,\omega_p)})\frac{\partial F_l}{\partial U_0^{(l)}}(\cdot,\omega_p),\label{dJmdu0} \\
\frac{\partial J_m}{\partial V^{\pm(l)}_n}= \sum_{p=1}^M (\overline{u(x_j,\omega_p)-F_l(\cdot,\omega_p)})\frac{\partial F_l}{\partial V^{\pm(l)}_n}(\cdot,\omega_p),\label{dJmdvn}
\end{align}
we denote here $F_j(U_0^{(j)},V^{\pm (j)}_1,\cdots,V^{\pm (j)}_{N_f},\omega_p,\kappa_1,\kappa_2,\kappa_3)$ by $F_j(\cdot,\omega_p)$ in order to simplify the notations.\\

Then, the algorithm follows the standard gradient method for $3+N_d (1+2N_f)$ variables. Once we have reconstructed the conductivity profile, i.e. the approximate values of $\kappa_1$, $\kappa_2$, $\kappa_3$, we can use (\ref{1mode}) again to calculate the approximate conductivity $\tilde{k}(\omega)$ by (\ref{empirical}) and the approximate $u_0$ by the following matrix formula, let $x\in \partial \Omega$,
\begin{equation}
\underbrace{\begin{pmatrix}
\tilde{u}(x,\omega_1) \\ 
\tilde{u}(x,\omega_2) \\ 
\vdots \\ 
\tilde{u}(x,\omega_M)
\end{pmatrix}}_{=\tilde{U}(x,\omega_1,\ldots,\omega_M)} \approx \underbrace{\begin{pmatrix}
q_0(\omega_1) & q(\widetilde{\lambda^+_1},\omega_1) & q(\widetilde{\lambda^-_1},\omega_1) & \cdots & q(\widetilde{\lambda^-_{N_f}},\omega_1)  \\ 
q_0(\omega_2) & q(\widetilde{\lambda^+_1},\omega_2) & q(\widetilde{\lambda^-_1},\omega_2) & \cdots & q(\widetilde{\lambda^-_{N_f}},\omega_2)  \\ 
\vdots & \vdots & \vdots & \ddots & \vdots  \\ 
q_0(\omega_M) & q(\widetilde{\lambda^+_1},\omega_M) & q(\widetilde{\lambda^-_1},\omega_M) & \cdots & q(\widetilde{\lambda^-_{N_f}},\omega_M) 
\end{pmatrix}}_{=L(\widetilde{\lambda^\pm_1},\ldots,\widetilde{\lambda^\pm_{N_f}},\omega_1,\ldots\omega_M)} \underbrace{\begin{pmatrix}
u_0(x) \\ 
v_1^+(x) \\ 
v_1^-(x) \\ 
\vdots \\ 
v_{N_f}^-(x)
\end{pmatrix}}_{=V(x)},
\end{equation}

%%%%%%%%%------------------------------

where $\tilde{u}(x,\omega)=u(x,\omega)-\frac{2}{\tilde{k}(\omega)+k_0}\mathfrak{f}(x)$, $q_0(\omega)=\frac{\tilde{k}(\omega)-k_0}{k_0(\tilde{k}(\omega)+k_0)}$, and $q(\widetilde{\lambda},\omega)=\frac{1}{k_0+\widetilde{\lambda}(\tilde{k}(\omega)-k_0)}-\frac{2}{\tilde{k}(\omega)+k_0}$. Then, the vector $V$ can be obtained by the formula
\begin{equation}\label{matriceL}
V(x) \approx (L^T L)^{\dag}L^T \tilde{U}(x,\omega_1,\cdots\omega_M),
\end{equation}
where $(L^T L)^{\dag}$ is the pseudo-inverse of the matrix $L^T L$. 
The conditioning of the matrix $L^T L$ depends in fact  on the distance between the sampling
values $\omega_j, \, j=1\,\cdots, M,$ and the frequency profile $\eqref{empirical}$.  
The approximate $u_0(x)$ is then recovered by taking the first coefficient of the vector $V(x)$.\\

Finally, the algorithm to reconstruct  $u_0$, can be summarized  in the following steps:

\begin{enumerate}
\item Give an apriori estimation $\widetilde{\lambda_n^\pm},\,  n= 1, \cdots, N_f$, of the 
eigenvalues $\lambda_n^\pm,\,  n= 1, \cdots, N_f$. 
\item Choose a step length $\alpha_m>0$ for the gradient descent.
\item Initialize the vectors $U_0|_0$, $V_1|_0, \cdots,V_n|_0$ and the coefficients $\kappa_1|_0$, $\kappa_2|_0$, $\kappa_3|_0$.
\item While $|\nabla J_m|$ is larger then a given threshold, we do
\begin{enumerate}
\item Calculate the values of the functions $F_j$ by (\ref{def_Fj}), and $\nabla J_m$ by (\ref{dJmdki}), (\ref{dJmdu0}), (\ref{dJmdvn}).
\item Update the parameters $\kappa_i|_{k+1}=\kappa_i|_k-\alpha_m \frac{\partial J_m}{\partial \kappa_i}$, $U_0^{(l)}|_{k+1}=U_0^{(l)}|_k-\alpha_m \frac{\partial J_m}{\partial U_0^{(l)}}$,  and $V_n^{\pm(l)}|_{k+1}=V_n^{\pm(l)}|_k-\alpha_m \frac{\partial J_m}{\partial V_n^{\pm(l)}}$.
\end{enumerate}
\item When $|\nabla J_m|$ is smaller then the threshold, we stop the iterations.
\item Use (\ref{matriceL}) with the approximate coefficients $\kappa_i$ obtained in the previous step to calculate the approximate value of $u_0(x)$ for every $x\in \partial \Omega$.
\end{enumerate}

%%%%%%%%%%%%%%%%%%%%%%%%%%%%%%%%%%%%%%%%%%%
\section{Reconstruction of the anomaly from $u_0$}
%%%%%%%%%%%%%%%%%%%%%%%%%%%%%%%%

In this section, we propose a numerical method to identify the  anomaly $D$
from  a finite number of Cauchy data of $(u_0(f_i), f_i), \, i= 1, \cdots, P, $ on $\partial \Omega$,
where $P\geq 1$.  We further assume that
the anomaly is located within an open subdomain $\Omega_0 \subset \Omega$ with $dist(\partial \Omega_0, \partial \Omega) \geq \delta_0 >0$ The scheme is based on the minimizing of a non convex functional
\[J(u)=\frac{1}{2} \int_{\partial \Omega}\sum_{i=1}^P |u-u_{meas}^{(i)}|^2 ds,\]
where $u_{meas}^{(i)}$ are the measured Dirichlet data corresponding to the $i$-th Neumann data and where $u$ is the solution to (\ref{maineq}) associated to the current domain $D \subset \Omega_0$. In our numerical simulations we take $P=2$  with $f_1=\langle e_1,\nu_{\Omega} \rangle$ and $f_2=\langle e_2, \nu_{\Omega} \rangle$, where $(e_1,e_2)$ is the canonical base of $\mathbb{R}^2$. \\

We further assume that  $D$ is within the class $\mathfrak D$, that is, it is star shaped and its boundary $\partial D$   can be described by the Fourier series:
\begin{equation}\label{modele dD}
\partial D=\left\{ X_0+ r(\theta)\begin{pmatrix}
\cos \theta \\ 
\sin \theta
\end{pmatrix} |\theta \in [0;2\pi) \right\}, \: r=\sum_{n=-N}^N c_n f_n,
\end{equation}

where $C=\begin{pmatrix}
c_{-N} \\ 
c_{-N+1} \\ 
\vdots \\ 
c_{N}
\end{pmatrix} \in \mathbb{R}^{2N+1} $, $f_n(\theta)=\cos(n\theta)$ for $0\leq n\leq N$ and $f_n(\theta)=\sin(n\theta)$ for $-N\leq n <0$. \\

Using (\ref{Pbuh}) in Theorem \ref{frechet},  and integration by parts, we have the expressions of the shape derivative corresponding to each Fourier coefficient $c_n$
\begin{equation}\label{dJexpression}
\frac{\partial J}{\partial c_n}=\int_{\Omega \setminus D} \nabla w \nabla u_h dX,
\end{equation}
for $-N\leq n \leq N$, 
where $h(\theta)=f_n(\theta) \langle \begin{pmatrix}
\cos \theta \\ 
\sin \theta
\end{pmatrix}  , \nu_D\rangle,$ and $w$ is the solution of the following equation
\begin{equation}
 \left\{
    \begin{array}{lr}
     \triangle w =0 & \text{in}\: \Omega \setminus \overline{D}, \\
     \dfrac{\partial w}{\partial \nu}=0 & \text{on}\: \partial D,\\
     \dfrac{\partial w}{\partial \nu}=u-u_{meas} & \text{on}\: \partial \Omega.
    \end{array}
  \right.
\end{equation}
Formula (\ref{dJexpression}) is also valid for the shape derivative corresponding to the displacement of $X_0$, in these cases, $h=\langle e_i, \nu_D \rangle$, $i=1,2$.\\
Those expressions are the basis of the following iterative algorithm:
\begin{enumerate}
\item Choose an initial domain $D_0$.
\item For each iteration, $i>0$:
      \begin{enumerate}
      \item Calculate the solution $u_i$  to (\ref{maineq}), associated to the domain $D_i$ for which the boundary $\partial D_i$ is calculated by (\ref{modele dD}).
      \item Calculate the shape derivatives $\frac{\partial J}{\partial x_1}$, $\frac{\partial J}{\partial x_2}$ and $\frac{\partial J}{\partial c_n}$ for all $-N\leq n\leq N$.
       \item Choose a step length $\alpha>0$ for the gradient descent.
      \item Update the parameters of the domain $X_{i+1}=X_i-\alpha \nabla_{X_0} J (X_i,C_i)$ and $C_{i+1}=C_i-\alpha \nabla_{C} J (X_i,C_i)$ with $\alpha>0$.
      \item If the updated domain is not entirely in $\Omega_0$ or if $R$ becomes negative, reduce the size of $\alpha$.
      \end{enumerate}
\item When $J(X_i,C_i)$ becomes smaller than a fixed threshold, we stop.
\end{enumerate}
%%%%%%%%%%%%%%%%%
\section{Numerical examples}
%%%%%%%%%%%%%%%%%
The numerical tests follow the steps presented here. All the  numerical experiments are done
using FreeFem++ \cite{FreeFem}.

\begin{figure}
    \centering
    \begin{subfigure}[b]{0.7\textwidth}
        \includegraphics[width=\textwidth]{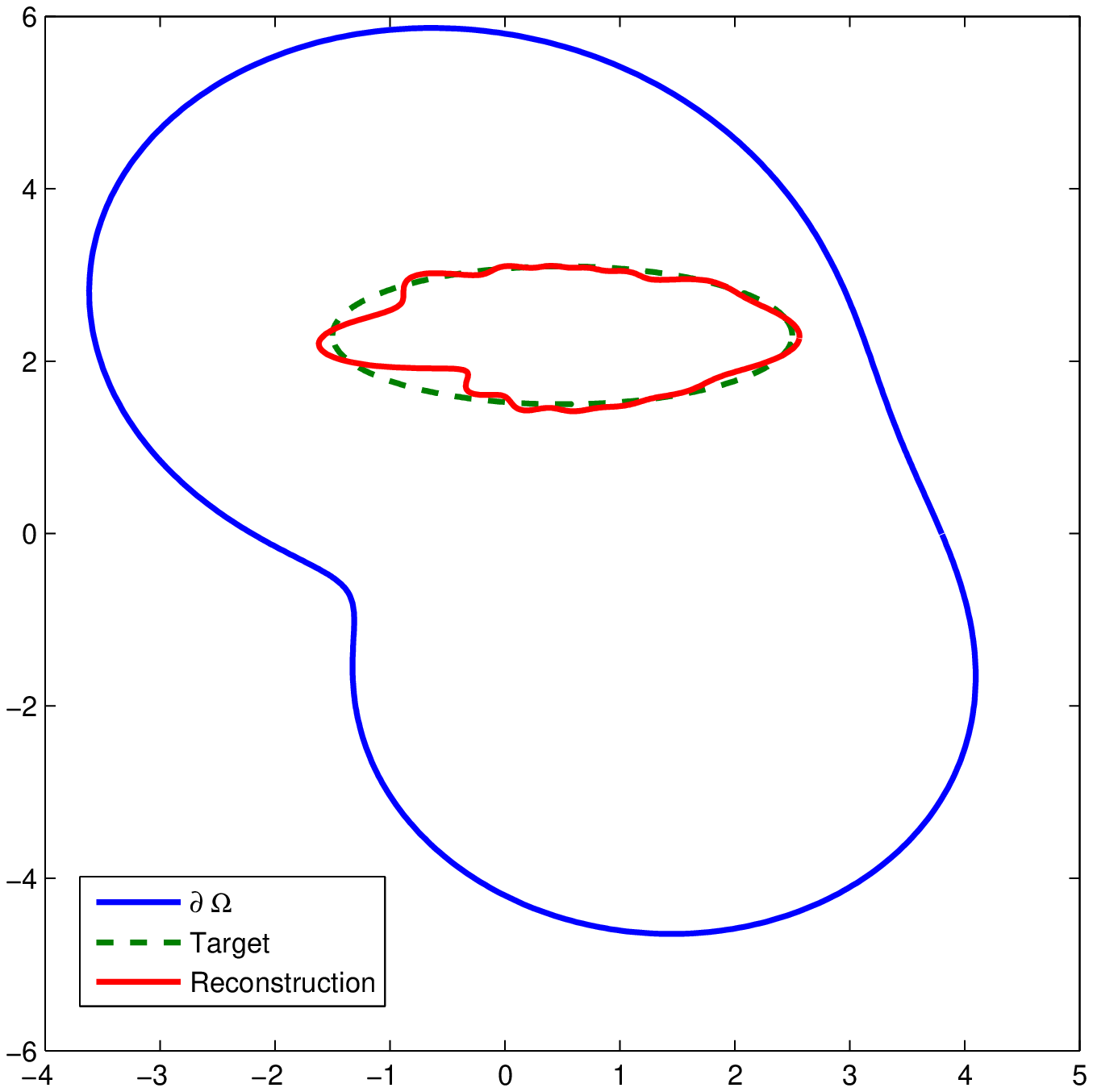}
        \caption{Target \& Reconstruction.}
        \label{figure_o2}
    \end{subfigure}
    \begin{subfigure}[b]{0.4\textwidth}
        \includegraphics[width=\textwidth]{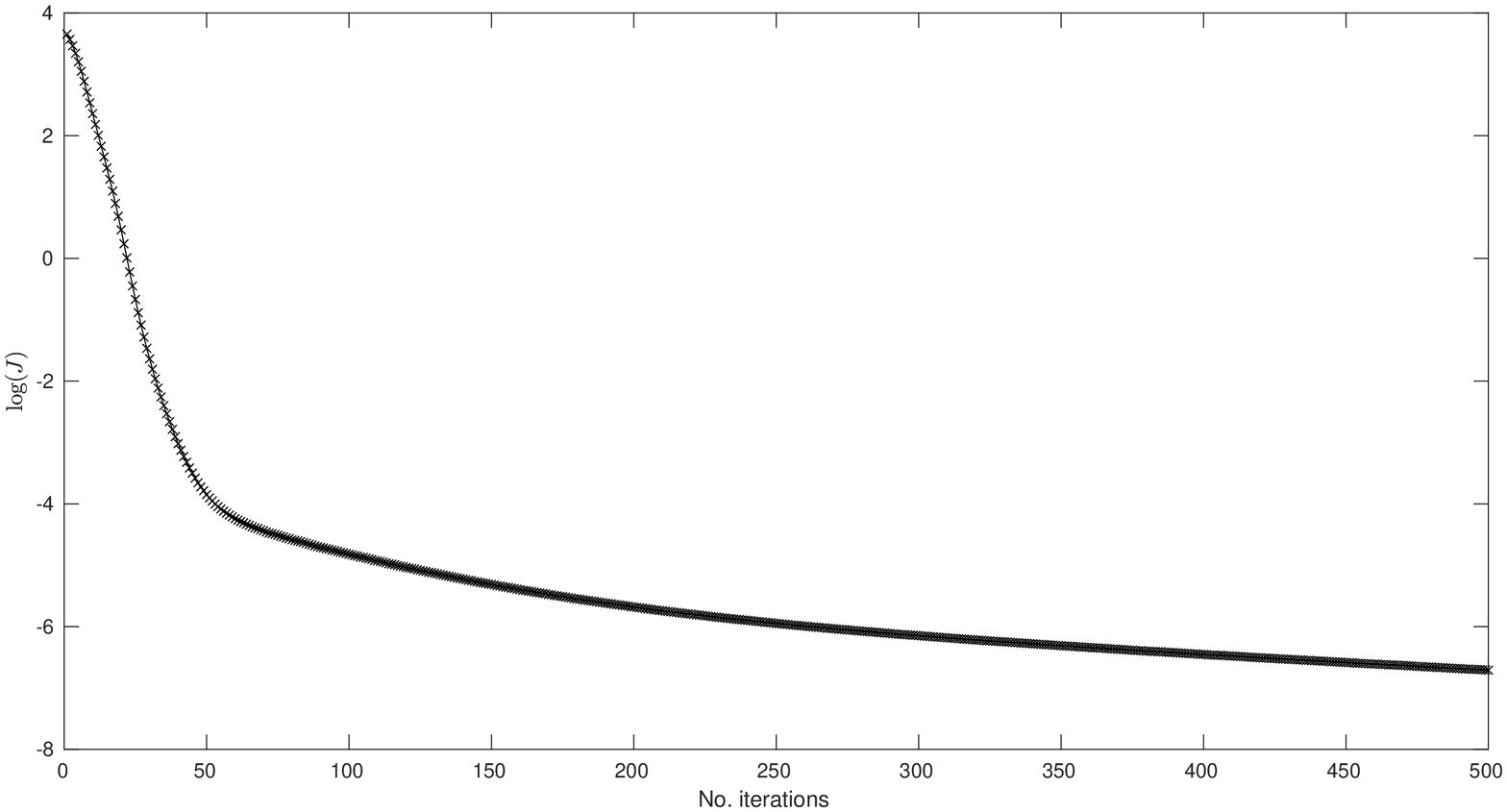} 
        \caption{$\log(J)$ with respect to Number of iterations.}
        \label{graphe_J_o2}
    \end{subfigure}
    \begin{subfigure}[b]{0.4\textwidth}
        \includegraphics[width=\textwidth]{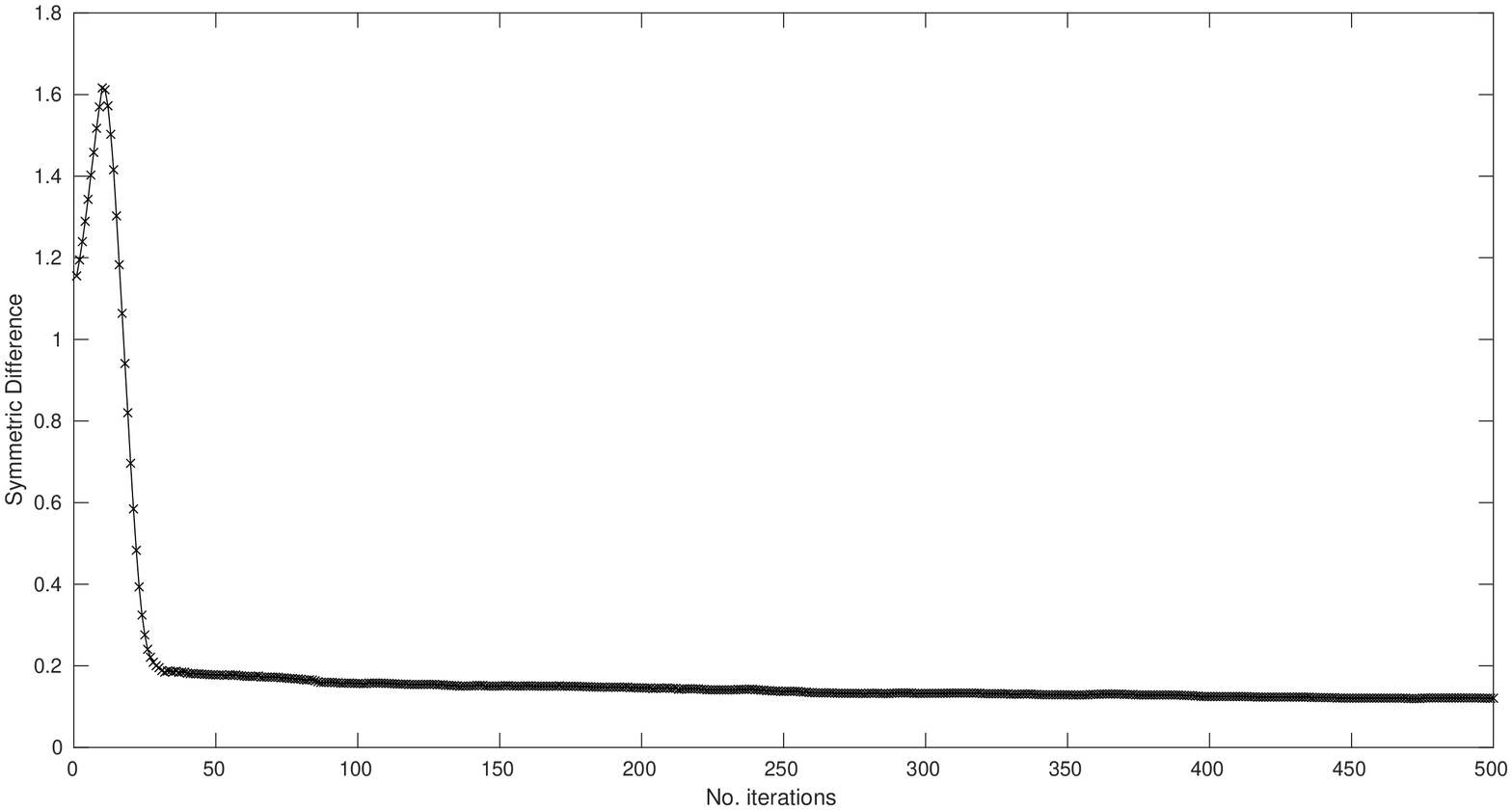} 
        \caption{Symmetric difference/Number of iterations.}
        \label{graphe_DS_o2}
    \end{subfigure}
    \begin{subfigure}[b]{0.4\textwidth}
        \includegraphics[width=\textwidth]{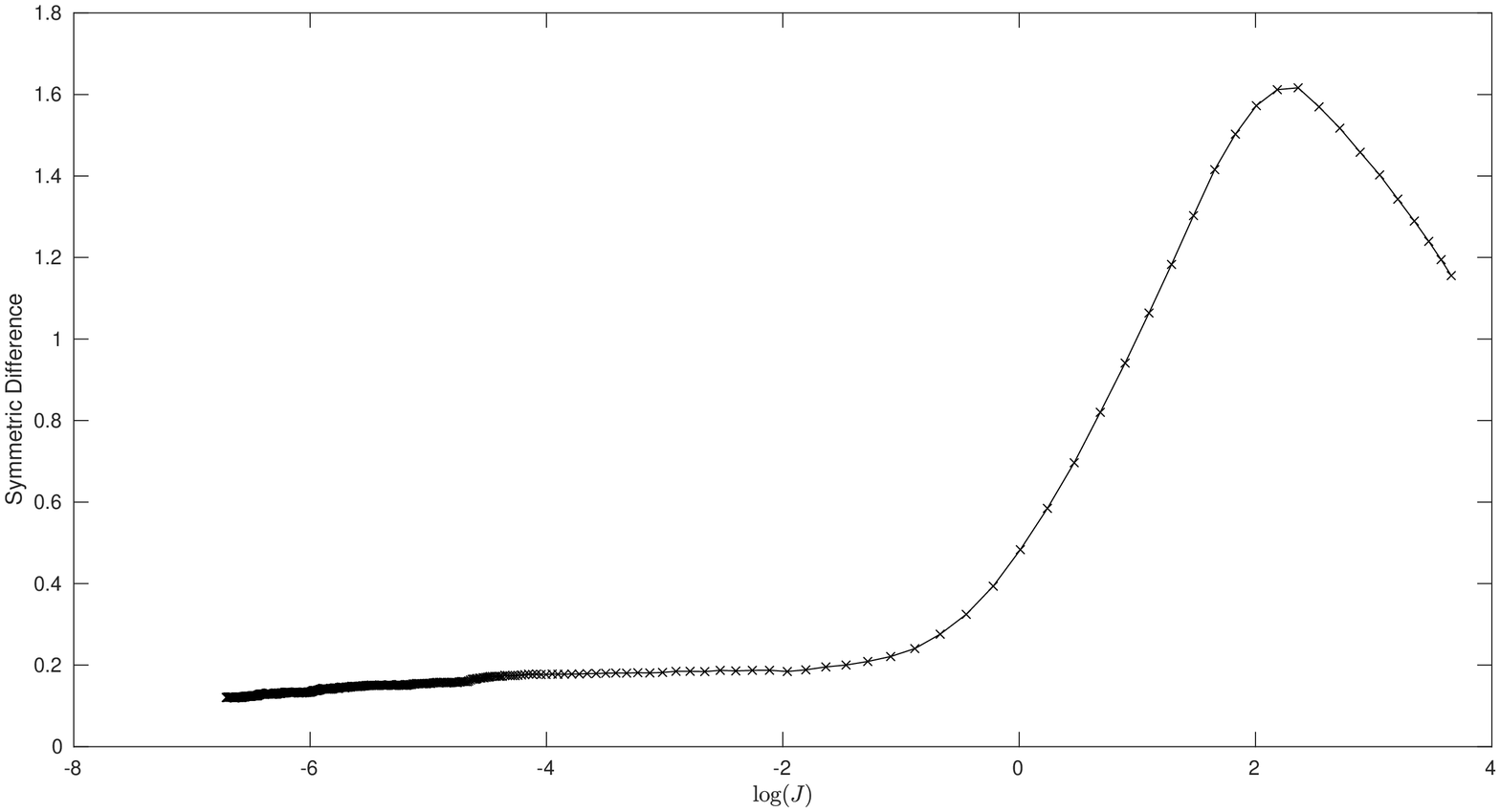} 
        \caption{Symmetric differences as a function of $\log(J)$.}
        \label{graphe_o2}
    \end{subfigure}
    \caption{Exemple 5: A non-ellptical shaped anomaly.}\label{o2}
\end{figure}

\begin{enumerate}
\item $\Omega$ is a centered ellipse defined by the equation: $\frac{x_1^2}{4^2}+\frac{x_2^2}{3^2} \leq 1$.
\item We use two linearly independent Neumann data: $f_1=\langle e_1,\nu_\Omega \rangle$ and $f_2=\langle e_2, \nu_\Omega\rangle$, where $(e_1,e_2)$ is the canonical base of $\mathbb{R}^2$.
\item The multifrequence conductivity follows the model (\ref{empirical}) with $\kappa_1=3$, $\kappa_2=2$, $\kappa_3=1$ and $\omega$ are integers from $1$ to $8$.
\item Only the first two eigenvalues are taken into consideration, and they are fixed as follows $\lambda^+_1=\frac{3}{4}$, $\lambda^-_1=\frac{1}{4}$ respectively in all cases.
\item In the algorithm to reconstruct $u_0$ and the conductivity profile, the initial guess of $u_0$ is the function $\mathfrak{f}$. 
\item The initial estimation of domain $D$ is a centered disk with a radius $\frac{1}{2}$.
\item We consider the first $15$ Fourier coefficients: $N=15$.
\item We use P1 finite elements for the numerical resolution of the PDEs.
\item At each iteration, we remesh the domain to adapt to the new predicted position and shape of the domain.
\item The algorithms stop if $J<10^{-5}$ or the number of iterations exceed $500$. All  the tests  have executed $500$ iterations.
\end{enumerate}
We present here several numerical simulations of the  proposed algorithm. We first give the errors in the reconstruction method of $u_0$  in Table (\ref{precision_u0}),
 and the errors in the reconstructed coefficients $\kappa_1$, $\kappa_2$, $\kappa_3$ in Table (\ref{kappas}). The errors are computed using the $L^2$-norm of the difference 
\bea
u_{0reconstruct}-u_0= error(u_{0reconstruct}):=\sqrt{\int_{\partial\Omega} |u_{0reconstruct}-u_0|^2 dx}.
\eea
 We show in the following figures the targets and the reconstruction result. We calculate also the relative symmetric difference $|D_i \bigtriangleup D_{target}|/|D_{target}|$ during the iterations, and we draw the curves of the symmetric difference with respect to $\log(u_i)$.  We finally give the relative symmetric difference of each shape  in Table \ref{DS_all}. Finally, we test  a reconstruction in domain
 $\Omega$  that has shape different from an ellipse  in Figure \eqref{o2}.

\begin{table}
\begin{center}
\begin{tabular}{|c|c|c|c|c|}
\hline 
• & ellipse & square & near-boundary & small-central \\ 
\hline 
$f=f_1$ & 0.04707 & 0.11973 & 0.00956 & 0.00502 \\ 
\hline 
$f=f_2$ & 0.01583 & 0.09905 & 0.02436 & 0.00893 \\ 
\hline 
\end{tabular} 
\end{center}
\caption{Errors between $u_{0reconstruct}$ and $u_0$.}\label{precision_u0}
\end{table}

\begin{table}
\begin{center}
\begin{tabular}{|c|c|c|c|c|c|}
\hline 
•&real value & ellipse & square & near-boundary & small-central \\ 
\hline 
$\kappa_1$&3 & 2.80971 & 3.36482 & 3.00287 & 6.65418 \\ 
\hline 
$\kappa_2$&2 & 1.79063 & 2.34197 & 1.96926 & 5.14671 \\ 
\hline 
$\kappa_3$&1 & 1.00212 & 0.987247 & 0.999658 & 1.13223 \\ 
\hline 
\end{tabular} 
\end{center}
\caption{Reconstructed constants in the frequency profile.}\label{kappas}
\end{table}

\begin{figure}
    \centering
    \begin{subfigure}[b]{0.4\textwidth}
        \includegraphics[width=\textwidth]{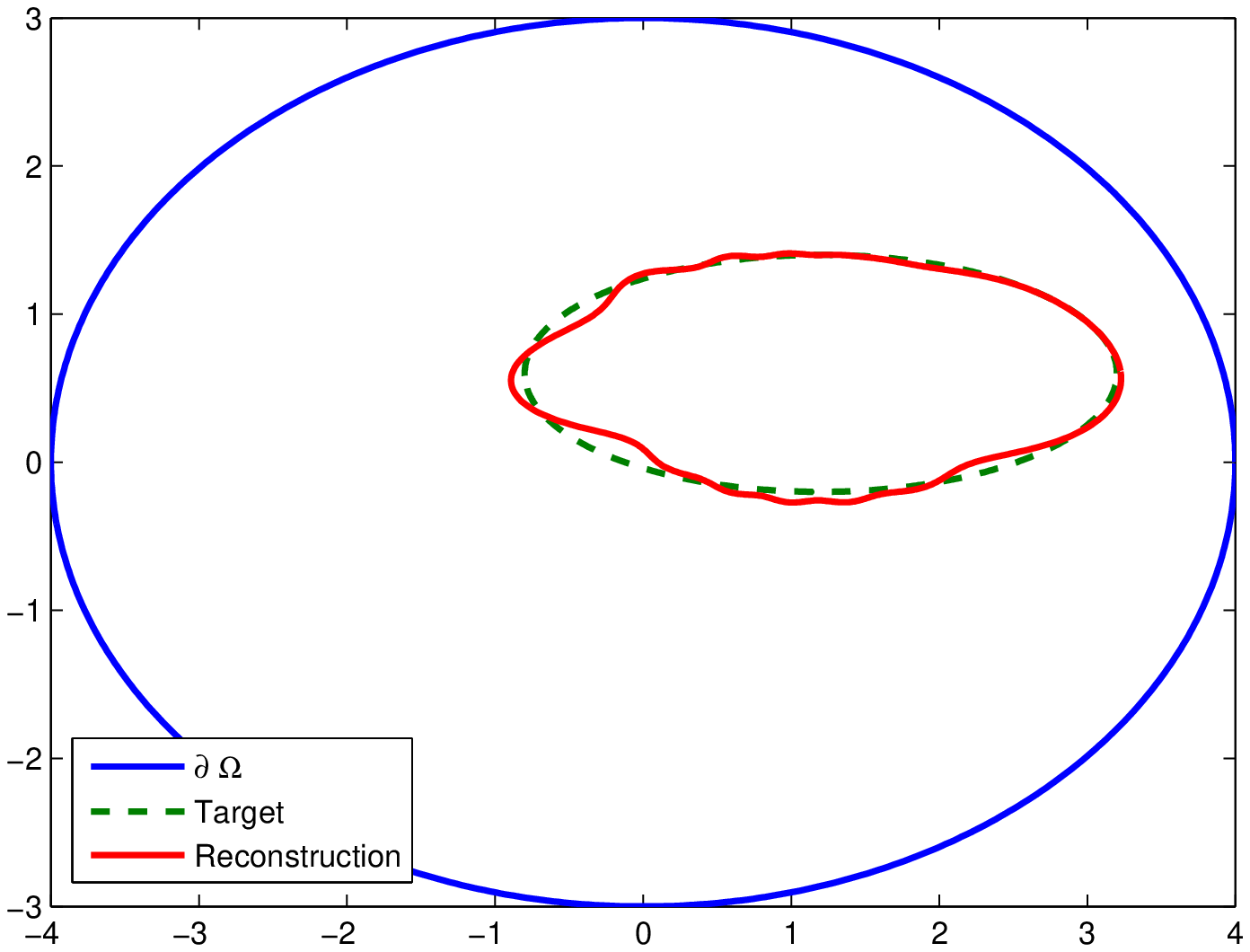}
        \caption{Target \& Reconstruction.}
        \label{figure_ellipse}
    \end{subfigure}
    \begin{subfigure}[b]{0.4\textwidth}
        \includegraphics[width=\textwidth]{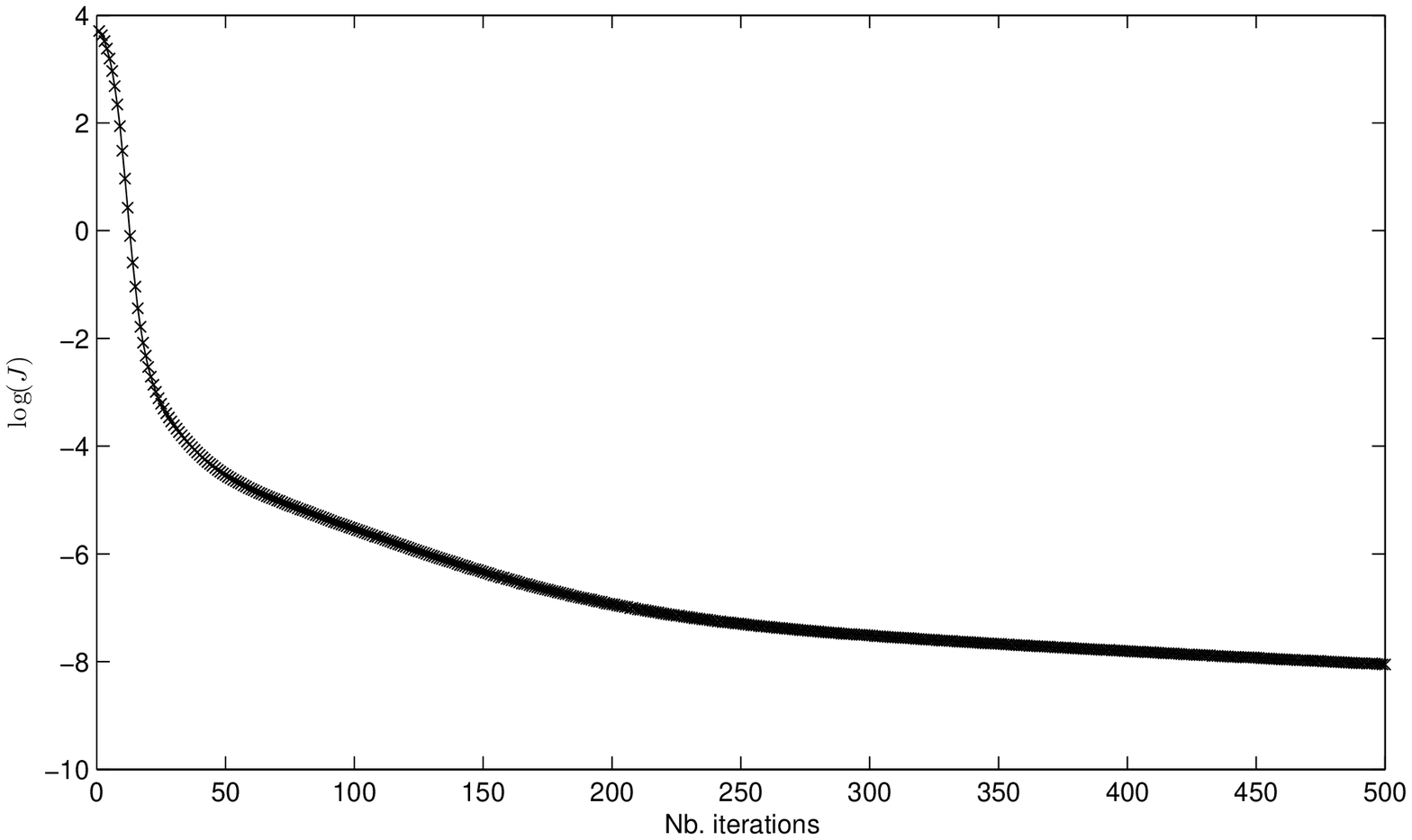} 
        \caption{$\log(J)$ during the algorithm.}
        \label{graphe_J_ellipse}
    \end{subfigure}
    \begin{subfigure}[b]{0.4\textwidth}
        \includegraphics[width=\textwidth]{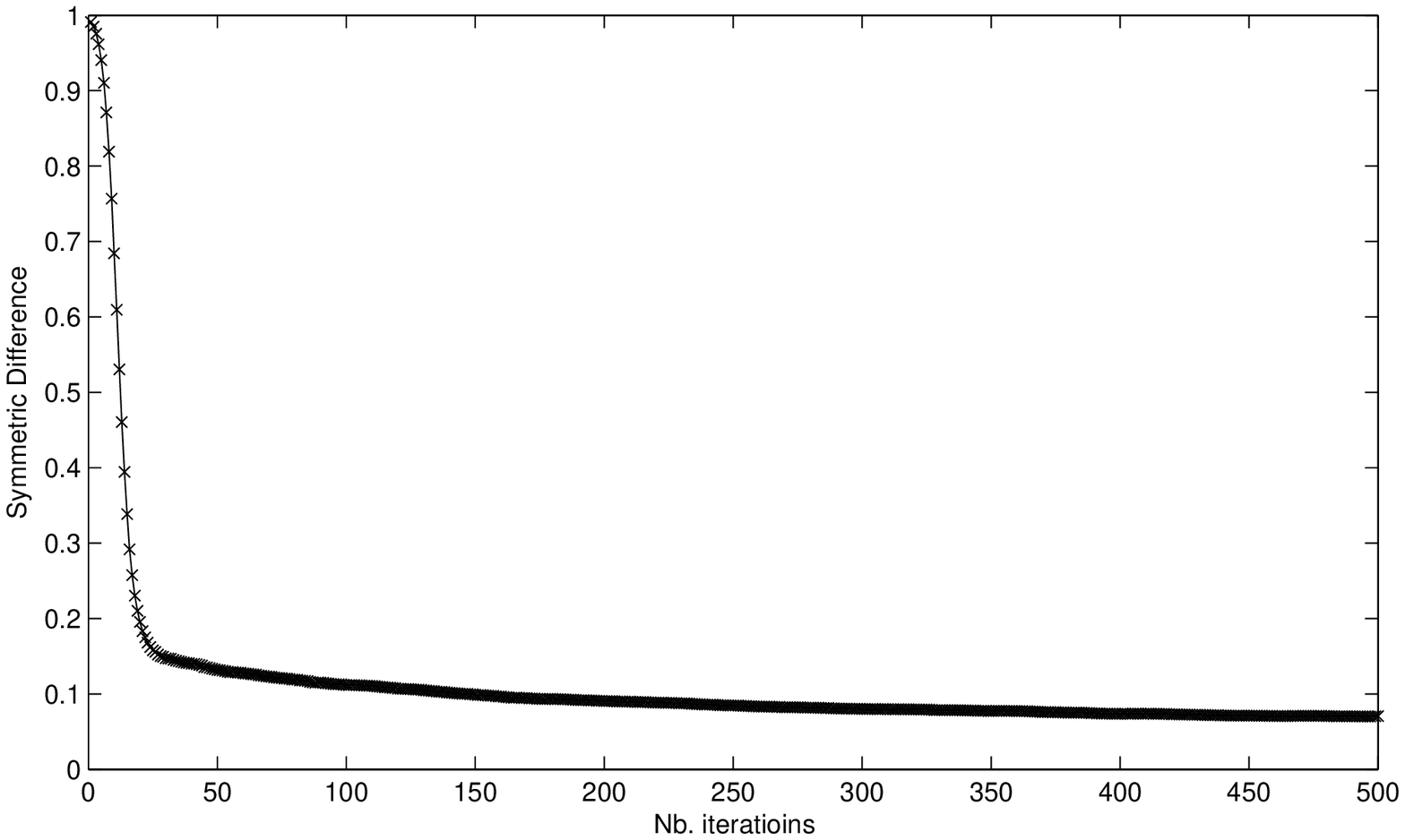} 
        \caption{Symmetric difference/Number of iteration.}
        \label{graphe_DS_ellipse}
    \end{subfigure}
    \begin{subfigure}[b]{0.4\textwidth}
        \includegraphics[width=\textwidth]{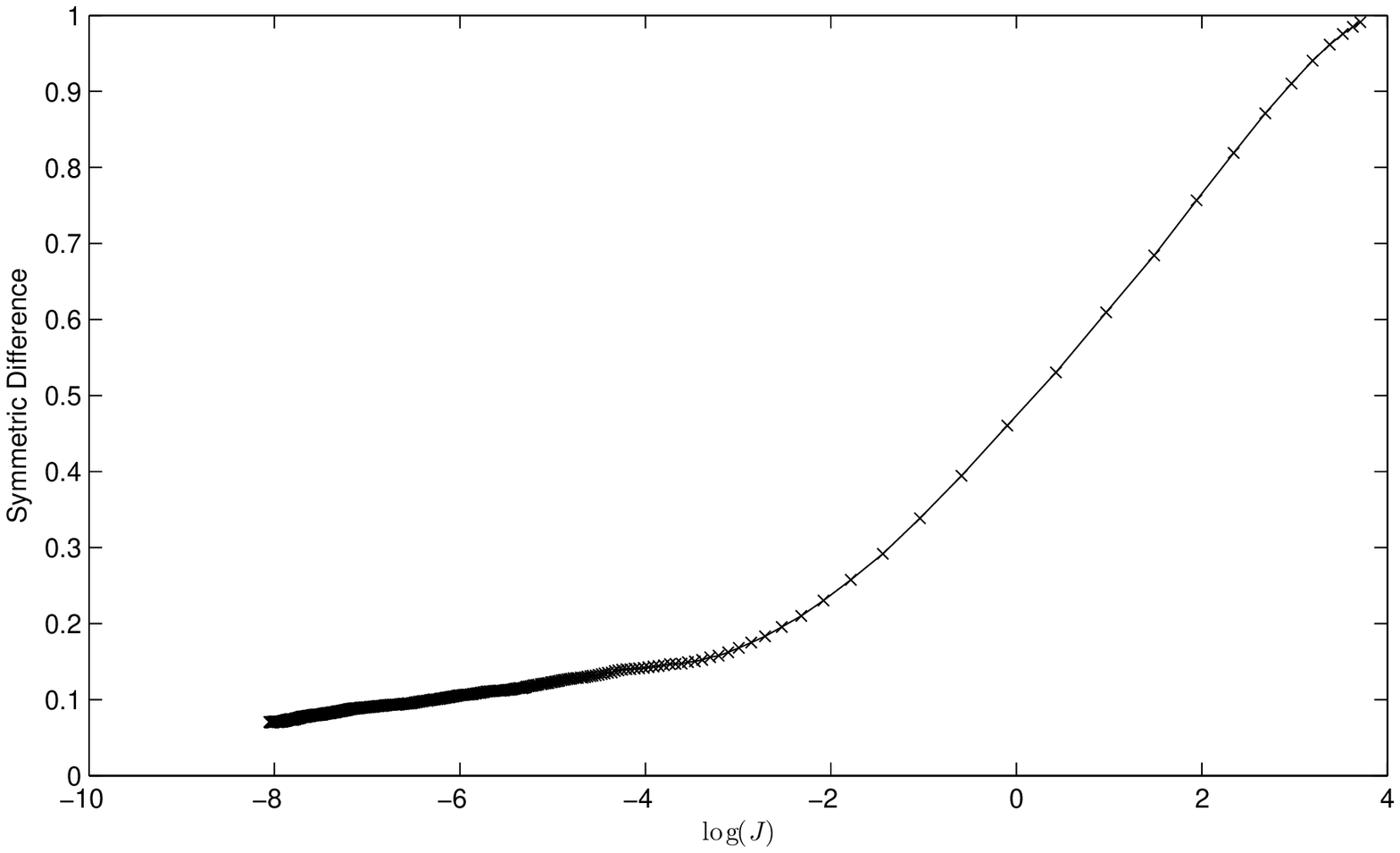} 
        \caption{Relation between $J$ and  symmetric differences.}
        \label{graphe_ellipse}
    \end{subfigure}
    \caption{Exemple 1: ellipse}\label{ellipse}
\end{figure}
\begin{figure}
    \centering
    \begin{subfigure}[b]{0.4\textwidth}
        \includegraphics[width=\textwidth]{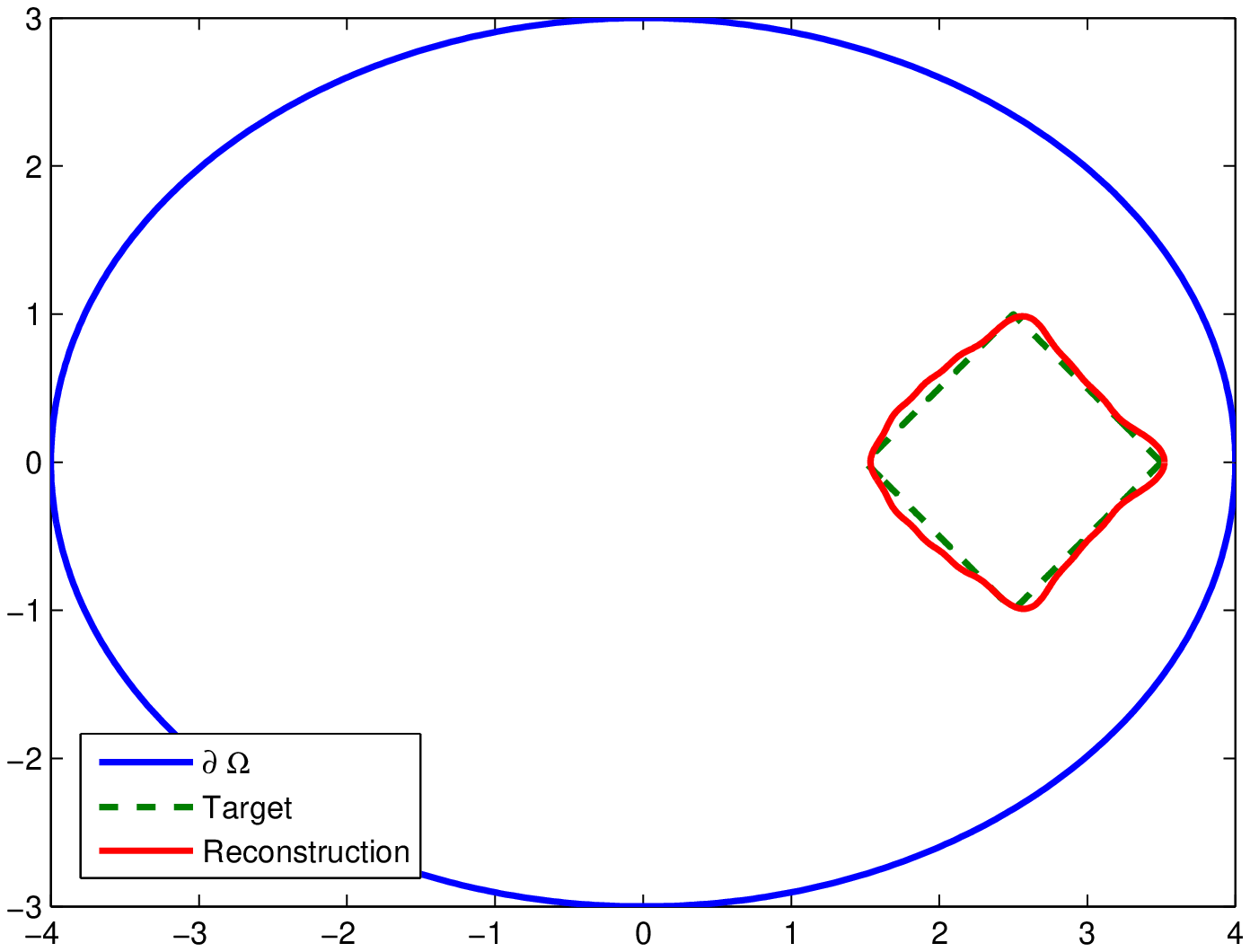}
        \caption{Target \& Reconstruction.}
        \label{target_carre}
    \end{subfigure}
    \begin{subfigure}[b]{0.4\textwidth}
        \includegraphics[width=\textwidth]{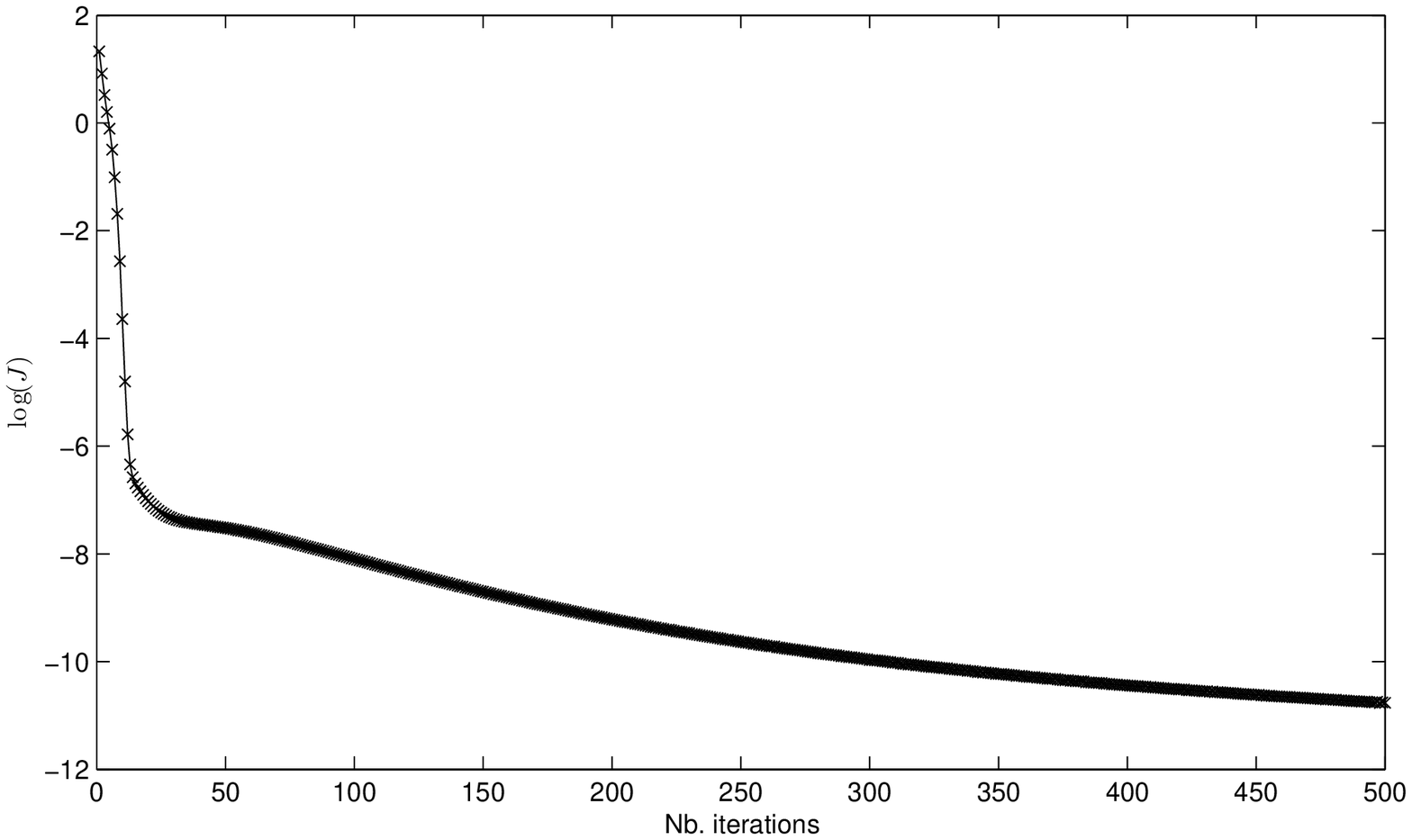} 
        \caption{$\log(J)$/Number of iterations. }
        \label{graphe_J_carre}
    \end{subfigure}
    \begin{subfigure}[b]{0.4\textwidth}
        \includegraphics[width=\textwidth]{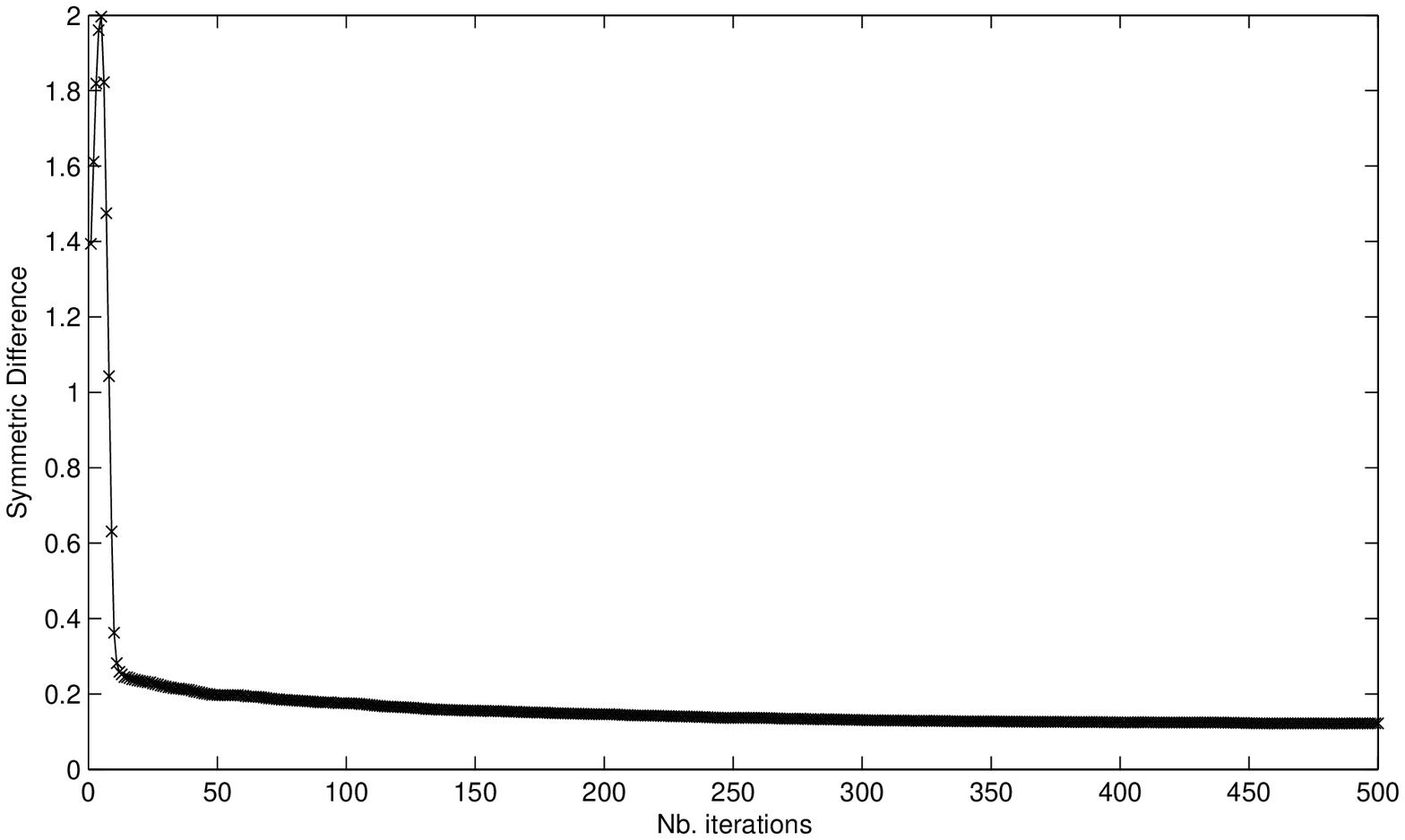} 
        \caption{Symmetric difference/Number of iterations.}
        \label{graphe_DS_carre}
    \end{subfigure}
    \begin{subfigure}[b]{0.4\textwidth}
        \includegraphics[width=\textwidth]{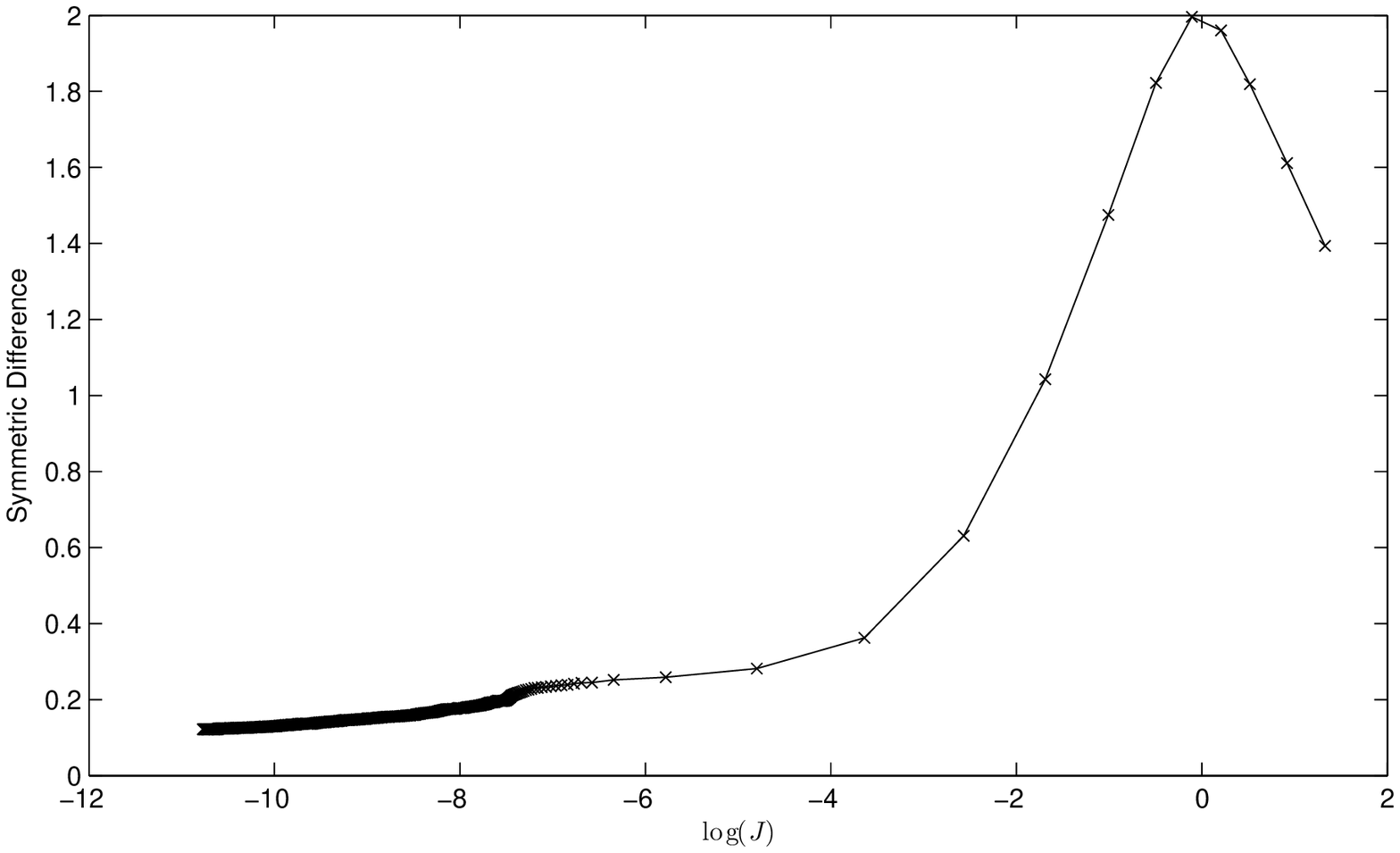} 
        \caption{Symmetric differences as a function of $\log(J)$.}
        \label{graphe_caree}
    \end{subfigure}
    \caption{Exemple 2: a square.}\label{caree}
\end{figure}
\begin{figure}
    \centering
    \begin{subfigure}[b]{0.4\textwidth}
        \includegraphics[width=\textwidth]{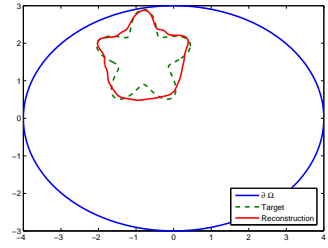}
        \caption{Target \& Reconstruction.}
        \label{target_hexag}
    \end{subfigure}
    \begin{subfigure}[b]{0.4\textwidth}
        \includegraphics[width=\textwidth]{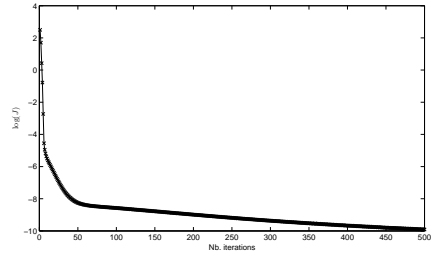} 
        \caption{$\log(J)$/Number of iterations.}
        \label{graphe_J_hexag}
    \end{subfigure}
    \begin{subfigure}[b]{0.4\textwidth}
        \includegraphics[width=\textwidth]{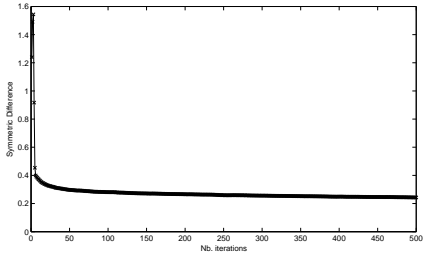} 
        \caption{Symmetric difference/Number of iterations.}
        \label{graphe_DS_hexag}
    \end{subfigure}
    \begin{subfigure}[b]{0.4\textwidth}
        \includegraphics[width=\textwidth]{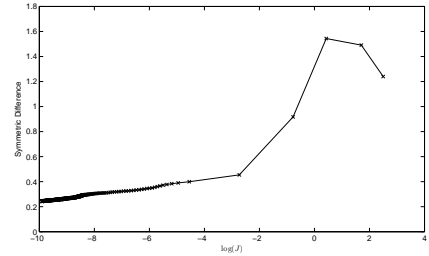} 
        \caption{Symmetric differences as a function of $\log(J)$.}
        \label{graphe_hexag}
    \end{subfigure}
    \caption{Exemple 3: a near boundary concave domain.}\label{hexag}
\end{figure}
\begin{figure}
    \centering
    \begin{subfigure}[b]{0.4\textwidth}
        \includegraphics[width=\textwidth]{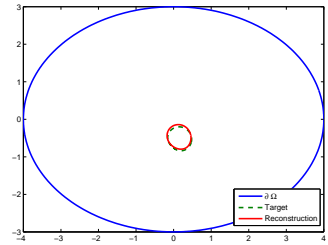}
        \caption{Target \& Reconstruction.}
        \label{figure_centre}
    \end{subfigure}
    \begin{subfigure}[b]{0.4\textwidth}
        \includegraphics[width=\textwidth]{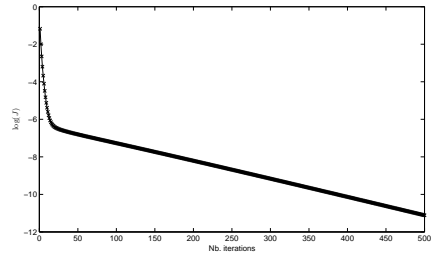} 
        \caption{$\log(J)$/Number of iterations.}
        \label{graphe_J_centre}
    \end{subfigure}
    \begin{subfigure}[b]{0.4\textwidth}
        \includegraphics[width=\textwidth]{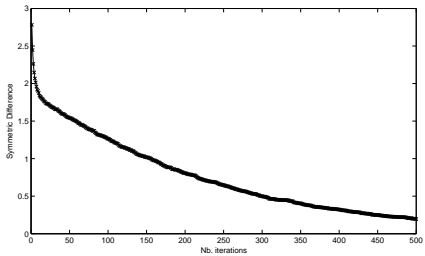} 
        \caption{Symmetric difference/Number of iterations.}
        \label{graphe_DS_centre}
    \end{subfigure}
    \begin{subfigure}[b]{0.4\textwidth}
        \includegraphics[width=\textwidth]{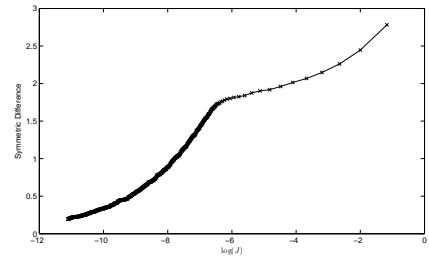} 
        \caption{Symmetric differences as a function of $\log(J)$.}
        \label{graphe_centre}
    \end{subfigure}
    \caption{Exemple 4: a centered small domain.}\label{centre}
\end{figure}
\begin{table}
\begin{center}
\begin{tabular}{|c|c|c|c|c|}
\hline 
• & ellipse & square & near-boundary & small-central \\ 
\hline 
$|D_i \bigtriangleup D_{target}|/|D_{target}|$ & 0.07055 & 0.12187 & 0.24299 & 0.19471 \\ 
\hline 
\end{tabular} 
\end{center}
\caption{Relative symmetric difference.}
\label{DS_all}
\end{table}
\section{Concluding remarks}

In this paper, by combining the spectral decomposition 
 derived in \cite{AT} and the linearization of the frequency independent part with respect to the
  shape of the anomaly, we have provided a new and efficient approach for reconstructing both the shape and conductivity parameter of a  conductivity anomaly from multifrequency boundary voltage measurements. The approach and results of this paper can be extended in several directions: (i) to reconstruct multiple anomalies from multifrequency boundary measurements; (ii) to investigate the reconstruction of anisotropic conductivity anomalies from multifrequency boundary  measurements, and (iii)  to study elastography imaging of visco-elastic anomalies.
  These new developments will be reported in forthcoming works. 

%%%%%%%%%%%%%%%%%%%%%%%%%%%%
\section{Acknowledgments}
This work has been partially supported by the 
LabEx PERSYVAL-Lab (ANR-11-LABX- 0025-01).

%%%%%%%%%%%%%%%%%%%%%%
%%%%%%%%%%%%%%%%%%%%%%

\end{document}